\documentclass[11pt]{amsart}
\usepackage{graphicx, color}
\usepackage{amssymb, amsmath, amsthm, mathdots}
\usepackage{youngtab}
\usepackage[all]{xy}
\usepackage{mathtools}
\usepackage{colonequals}
\begin{document}
\newtheorem{defn0}{Definition}[section]
\newtheorem{prop0}[defn0]{Proposition}
\newtheorem{thm0}[defn0]{Theorem}
\newtheorem{lemma0}[defn0]{Lemma}
\newtheorem{coro0}[defn0]{Corollary}
\newtheorem{exa}[defn0]{Example}

\numberwithin{equation}{section}
\def\rig#1{\smash{ \mathop{\longrightarrow}
    \limits^{#1}}}
\def\nwar#1{\nwarrow
   \rlap{$\vcenter{\hbox{$\scriptstyle#1$}}$}}
\def\near#1{\nearrow
   \rlap{$\vcenter{\hbox{$\scriptstyle#1$}}$}}
\def\sear#1{\searrow
   \rlap{$\vcenter{\hbox{$\scriptstyle#1$}}$}}
\def\swar#1{\swarrow
   \rlap{$\vcenter{\hbox{$\scriptstyle#1$}}$}}
\def\dow#1{\Big\downarrow
   \rlap{$\vcenter{\hbox{$\scriptstyle#1$}}$}}
\def\up#1{\Big\uparrow
   \rlap{$\vcenter{\hbox{$\scriptstyle#1$}}$}}
\def\lef#1{\smash{ \mathop{\longleftarrow}
    \limits^{#1}}}
\def\O{{\mathcal O}}
\def\L{{\mathcal L}}
\def\M{{\cal M}}
\def\N{{\cal N}}
\def\P#1{{\mathbb P}^#1}
\def\PP{{\mathbb P}}
\def\CC{{\mathbb C}}
\def\QQ{{\mathbb Q}}
\def\Z{{\mathbb Z}}
\def\H{{\mathcal H}}
\def\sym{{\mathrm{Sym}}}
\newcommand{\defref}[1]{Def.~\ref{#1}}
\newcommand{\propref}[1]{Prop.~\ref{#1}}
\newcommand{\thmref}[1]{Thm.~\ref{#1}}
\newcommand{\lemref}[1]{Lemma~\ref{#1}}
\newcommand{\corref}[1]{Cor.~\ref{#1}}
\newcommand{\exref}[1]{Example~\ref{#1}}
\newcommand{\secref}[1]{Section~\ref{#1}}
\newcommand{\magma}{\textsf{Magma}}
\def\geq{\geqslant}
\def\leq{\leqslant}
\def\ge{\geqslant}
\def\le{\leqslant}

\newcommand{\qedd}{\hfill\framebox[2mm]{\ }\medskip}
\newcommand{\codim}{\textrm{codim}}
\newtheorem{rem0}[defn0]{Remark}
\def\proof{{\it Proof:\  \ }}
\author[C. Ciliberto, G. Ottaviani \MakeLowercase{with an appendix by} J. Caro, J. Duque-Rosero]{Ciro Ciliberto, Giorgio Ottaviani\smallskip\\
\MakeLowercase{with an appendix by} Jerson Caro, Juanita Duque-Rosero}

\address{Ciro Ciliberto\\
Dipartimento di Matematica, Universit\`a di Roma Tor Vergata, Via della Ricerca Scientifica, 00173 Roma, Italy}
\email{cilibert@mat.uniroma2.it}

\address{Giorgio Ottaviani\\
 Dipartimento di Matematica e Informatica ``Ulisse Dini'', Universit\`a di Firenze,  Viale Morgagni, 67/A, 50134 Firenze, Italy}
\email{giorgio.ottaviani@unifi.it}

\address{Jerson Caro\\
Department of Mathematics and Statistics, Boston University, 665 Commonwealth Ave, Boston, MA 02215, USA}
\email{jlcaro@bu.edu}

\address{Juanita Duque-Rosero\\
Department of Mathematics and Statistics, Boston University, 665 Commonwealth Ave, Boston, MA 02215, USA}
\email{juanita@bu.edu}

\title{The general ternary form can be recovered by its Hessian}
\maketitle

\begin{abstract}
The Hessian map is the rational map that  sends a homogeneous polynomial to the determinant of its Hessian matrix.
We prove that the Hessian map is birational on its image for ternary forms of degree $d\ge 4$, $d\neq 5$,
by considering the action of  the orthogonal group. In a previous paper we proved the analogous result for binary forms,
with more geometric techniques.

\end{abstract}

\section{Introduction}

In \cite{CO} we studied the Hessian map, which is the rational map

 $$h_{d,r}\colon\PP(\sym^d(\CC^{r+1}))\dashrightarrow\PP(\sym^{(d-2)(r+1)}(\CC^{r+1}))$$ 
sending a form  $f$ to its Hessian ${\rm Hess}(f)$,  which is the determinant of the Hessian matrix. We proved that $h_{d,1}$ is birational on its image for any $d\ge 5$ (which is sharp)
and that $h_{3,3}$ is birational on its image, while it is classical that $h_{3,2}$ is generically $3\colon 1$.
In this paper we prove:

\begin{thm0}\label{thm:ternary_even} The Hessian map $h_{d,2}$ is birational on its image for any $d\ge 4$, $d\neq 5$.
\end{thm0}

We think the assumption $d\neq 5$ can be removed but our proof does not work in this case. The proof considers separately the cases $d$ even and odd.

Our technique uses group actions. The Hessian map is equivariant for the action of ${\rm SL}(r+1)$. This action  has a unique closed orbit, both in 
the source and in the target space, which is the Veronese variety. This does not help too much since the Hessian map is not defined on the Veronese variety.
Our strategy consists in considering the action of a smaller group, namely the orthogonal group ${\rm SO}(r+1)$ with respect to a nondegenerate  quadratic form $q$.
This action has finitely many closed orbits, according to the fundamental harmonic decomposition, and when the degree $d$ is even has a unique fixed point, corresponding
to the form $q^{d/2}$. The Hessian map is defined at this point, it sends the unique fixed point  in  the  domain  to the unique fixed point in the target and it has maximal rank at this point in most cases, precisely when a quadratic equation, made explicit in Proposition \ref{prop:maxrank}, has no integer solutions. There exist certain values of $d$
when there are integer solutions of this equation (the first one for $r=2$ is $d=14$). In these cases our analysis is more subtle, and we need to analyze the Hessian map at a second point,
which is  of the form $q^{d/2-1}\ell^2$, where $\ell$ is an isotropic linear form, namely  the hyperplane $\ell=0$ is tangent to the quadric $Z(q)$, the zero locus of $q$. 
The rank of the Hessian map at this second point is maximal, at least for $r=2$ and $d\ge 6$, since the candidate values to fail the maximality are solutions of a cubic equation, made explicit in Proposition \ref{prop:maxrankql2},  which has no integer solutions in this range thanks to Theorem \ref{thm:integralPts2} of the Appendix.
These consideration allow to prove the birationality after an analysis of the resolution of the indeterminacy of the Hessian map, that we are able to do only for ternary forms. 
In the case $d=2k+1$ odd we perform a similar analysis around the form $q^k\ell$ where $\ell$ is an isotropic linear form.

We outline the content of the paper. In the second section we study the Hessian map around the closed orbits of the ${\rm SO}(r+1)$-action on $\PP(\sym^d(\CC^{r+1}))$.
We state Theorem \ref{thm:criterion}, which is a criterion for the birationality of the Hessian map based on which closed orbits belong
to the graph  of the Hessian map.
In the third section we analyze the indeterminacy values of the Hessian map for ternary forms, approaching the powers $\ell^d$ with $\ell$ an  isotropic linear form, 
which make the only closed  ${\rm SO}(3)$-orbit   where the Hessian map is not defined. In Theorem \ref{thm:main} we prove our main Theorem \ref{thm:ternary_even}.

The appendix by Jerson Caro and  Juanita Duque-Rosero shows that some cubic equations have no integer solutions beyond a few ones. This allows to show
that the differential of the Hessian map has maximum rank at some relevant points. This is a crucial step for our technique. Without this appendix
we could prove only that Theorem 1.1 holds with at most finitely many exceptions for $d$.

For basic facts on the Hessian see \cite{Ru}, for an interesting recent approach see also \cite{BFP}. 

{ After this paper has been written we received the preprint \cite{Beo} by V. Beorchia, where some results related to our Theorem 1.1 are proved.}

{ GO and JD-R thank Frank Sottile for the organization of the Texas Algebraic Geometry Symposium at Texas A\&M in April 2024, where discussions began that led to the appendix taking shape. }
{ CC and GO are members of Italian GNSAGA-INDAM. GO's research is funded by the European Union - Next Generation EU, M4C1, CUP B53D23009560006, PRIN 2022- 2022ZRRL4C - Multilinear Algebraic Geometry. JC was supported by Simons Foundation grant no. 550023 and JDR was partially supported by Simons Foundation grant no. 550023.}
\section {The Hessian map and the ${\rm SO}(r+1)$-action on $\PP(\sym^d(\CC^{r+1}))$}
 Set $V=\CC^{r+1}$ with  coordinates $(x_0,\ldots, x_r)$ and 
$${\rm Hess}(f)\colon = \det \left(\frac{\partial^2 f}{\partial x_i\partial x_j}\right)_{0\leq i,j\leq r}.$$
The Hessian map defined in the Introduction is ${\rm SL}(V)$-equivariant, this means that if $g\in {\rm SL}(V)$ and $f\in{\rm Sym}^d(V)$ then 
$${\rm Hess}(g\cdot f)=g\cdot {\rm Hess}(f),$$
where the actions of ${\rm SL}(V)$ on forms is the obvious one. 

Any nondegenerate quadratic form $q\in {\rm Sym}^2(V)$ defines the special orthogonal group
${\rm SO}(V,q)=\{g\in {\rm SL}(V)| g\cdot q=q\}$. In particular the Hessian map is also ${\rm SO}(V,q)$-equivariant for any nondegenerate quadratic form $q$.
Under the ${\rm SO}(V,q)$-action the space ${\rm Sym}^d(V)$ is not irreducible and splits as \cite[Exerc. 19.21]{FH}, \cite[Cor. 5.2.5]{GW}
\begin{equation}\label{eq:harmdec} \sym^d(V) = \H_d\oplus q\H_{d-2}\oplus q^2\H_{d-4}\oplus\ldots,\end{equation} where $\H_i$ is the space of  \emph{harmonic homogeneous polynomials of degree $i$}, that by definition are killed by the dual operator $q^*$.
 If $q=\sum_{i=0}^r x_i^2$, the operator $q^*$ is the classical Laplacian operator $\Delta=\sum_{i=0}^r\frac{\partial^2}{\partial x_i^2}$. The decomposition (\ref{eq:harmdec}) is known as the {\it harmonic decomposition}.
Note that the power of a linear form $\ell^d$ is harmonic if and only if  the linear form  $\ell$ is isotropic, namely the hyperplane $\ell=0$  is tangent  to the isotropic quadric  $Z(q)=\{q=0\}$. 
 For a geometric point of view see section 1.4.5 in \cite{Dolg}.  We summarize the main properties of the harmonic decomposition in the following well known Proposition,
{ which we state separately from (\ref{eq:harmdec}) just to give a few definitions that will be used in the rest of the paper and to show in the proof an easy algorithm to construct the harmonic decomposition.}
\begin{prop0}\label{prop:harmdec}
Fix a nondegenerate quadratic form  $q\in \sym^2(V)$  and its dual  $q^*\in\sym^2 (V^\vee)$.
For any $f\in\sym^d(V)$ there exist unique  $f_{d-2i}\in\sym^{d-2i}(V)$   for $i=0,\ldots, \lfloor d/2\rfloor$ such that
$$f=\sum_{i=0}^{\lfloor d/2\rfloor}q^if_{d-2i},\qquad q^*(f_{d-2i})=0,\quad \forall i=0,\ldots ,\lfloor d/2\rfloor $$
The polynomial $f_d$ is called the  \emph{harmonic part}  of $f$.
The polynomial $f_j$ is called the   \emph{$j$-th harmonic summand} of $f$.
The harmonic decomposition of $q^jf$ is obtained by multiplying by $q^j$ the harmonic decomposition of $f$.
\end{prop0}
\begin{proof}
The harmonic part $f_d$ can be found as $f_d=f-qg$ where $g$ is the unique solution of the linear system  $q^*(f)=q^*(qg)$. 
The following summands can be found recursively, starting from $f_{d-2}$ which is the harmonic part of $g$.\qedd
\end{proof}
\vskip 0.4cm 

\begin{rem0}\label{rem:bw} {\rm  The bilinear symmetric form $q$ on $V$ can be extended to a bilinear symmetric form $Q$ on $\sym^d(V)$   by the condition $Q(v^d,w^d)=\left(q(v,w)\right)^d$,
$\forall v, w\in V$,
which is called the  \emph{Bombieri-Weyl form}.  The summands of the harmonic decomposition (\ref{eq:harmdec}) are orthogonal with respect to $Q$. The linear map that sends a polynomial to its $j$-th harmonic summand coincides with the orthogonal projection on the summand $\H_j$, in particular it is $SO(V)$-equivariant.}
\end{rem0}

Our arguments do not depend on the choice of $q$, so we denote often ${\rm SO}(V)={\rm SO}(V,q)$, the form $q$ being understood. 
The closed ${\rm SO}(V)$-orbits in  $\PP(\sym^d (V))$  are finitely many, precisely there is one closed orbit in each summand of (\ref{eq:harmdec}).
The closed orbit in the summand $q^i\H_{d-2i}$ consists of  polynomials $q^i\ell^{d-2i}$ where $\ell$ is any linear form such that the hyperplane $\ell=0$ is tangent to  $Z(q)$. 
They are isomorphic to   the dual of  $Z(q)$ (isomorphic in turn to $Z(q)$)  embedded by the $(d-2i)$-Veronese embedding, unless $d$ is even and $d=2i$, when the orbit consists of the single point
$q^{d/2}$. Let us denote by $\mathbb Q_{r-1}$ a smooth quadric in $\PP(V)$.  Summing up we have:

\begin{lemma0}\label{lem:isoquad}
 The polarized quadric $(\mathbb Q_{r-1},\O(d))$ is the unique closed ${\rm SO}(V)$-orbit  in the space $\PP (\H_d)$. There is a natural isomorphism of  ${\rm SO}(V)$-modules 
$$H^0(\mathbb Q_{r-1},\O(d)) = \H_d.$$
There is a unique one-dimensional summand $\langle q^k\rangle$ in (\ref{eq:harmdec}) exactly when $d=2k$ is even. 
\end{lemma0}

\begin{prop0}\label{prop:hescon} Let $\dim (V)=r+1$, $d=2k$. For any nondegenerate quadratic form $q$ we have 
$${\rm Hess}(q^k)=cq^{(r+1)(k-1)}$$ for some nonzero scalar $c$, which will be computed  in Proposition \ref{prop:ql}.
\end{prop0}

\begin{proof} The polynomial ${\rm Hess}(q^k)$  must be invariant by the ${\rm SO}(V,q)$-action, hence it is a power of $q$, up to a scalar multiple.
The fact that $c\neq 0$ can be checked with the diagonal form $q=\sum_{i=0}^rx_i^2$. Indeed,  the monomial power $x_0^{2(r+1)(k-1)}$ appears
as a summand in the Hessian and the only entries of the Hessian matrix that contain a monomial power of $x_0$ are the diagonal ones, that contain $2kx_0^{2(k-1)}$,  which are all nonzero. \qedd  \end{proof}
\vskip 0.4cm

\begin{rem0} {\rm
The argument in the proof of Proposition \ref{prop:hescon} is well known and may be applied all the times
we have an action of a group $G$ with a single generator $F$ of degree $d$ of the invariant ring $\left(\oplus_{i}{\rm Sym}^i(V)\right)^G$.
 Then the Hessian of $F$ is equal to $cF^{(r+1)(d-2)/d}$ for a scalar $c$ 
when $\frac{(r+1)(d-2)}{d}\in\Z$ and it is zero when $\frac{(r+1)(d-2)}{d}\notin\Z$. 
In these cases the degree of the polar map (see \cite{CRS}) is one if $c\neq 0$  ($f$ is homaloidal)
or it is zero if $c=0$ ($f$ has vanishing Hessian). This case has been called {\it totally Hessian} in \cite[Rem. 3.5]{CRS}.
Note that $V$ is prehomogenous for the action of $\CC^*\times {\rm SO}(V)$, see \cite[Example 19]{Man}. 

 The simplest example of this
behaviour is when $F$ is the  symmetric $n\times n$ determinant, then ${\rm Hess}(F)=cF^{(n+1)(n-2)/2}$ for a nonzero scalar $c$.  This was noticed by Beniamino Segre in \cite{Seg}.

When $F$ is the determinant of a matrix in $n^2$ indeterminates, its Hessian is
$cF^{n(n-2)}$ for a nonzero scalar $c$. In the same way ${\rm Hess}(F^k)=cF^{n(nk-2)}$ for another nonzero scalar $c$.
When $F$ is the $n\times n$ Pfaffian, then ${\rm Hess}(F)=cF^{(n-1)(n-4)/2}$.}

\end{rem0}

\begin{prop0}\label{prop:isotropy} Let $\mathbb Q=\mathbb Q_{r-1}\subset\PP(V)$ be the isotropic quadric.  Fix $x\in \mathbb Q$ and let $P_x=\{g\in {\rm SO}(V)| gx=x\}$. Let $H=x^\perp$ be the tangent hyperplane to $\mathbb Q$ at $x$.

(i) The action of ${\rm SO}(V)$ on $\PP(V)$ has exactly two orbits, namely $\PP(V)\setminus \mathbb Q$ and $\mathbb Q$.

(ii)
The action of $P_x$ on $\PP(V)$ has exactly five orbits for $r\ge 4$, namely $\PP(V)\setminus \left( \mathbb Q\cup H\right)$, $\mathbb Q\setminus H$, $H\setminus \mathbb Q$,  $(\mathbb Q\cap H)\setminus \{x\}$ and $\{x\}$. 

For $r=3$ the orbits are six, since $\mathbb Q\cap H$ consists of two lines.

For $r=2$ the orbits are four since  $(\mathbb Q\cap H)\setminus \{x\}=\emptyset$.

For $r=1$ the orbits are three since $H\setminus \mathbb Q = (\mathbb Q\cap H)\setminus \{x\} = \emptyset$.

\end{prop0}
\begin{proof} It is well known and a straightforward computation in a coordinate system. \qedd
\end{proof}\medskip

Proposition \ref{prop:hescon} has the following generalization to the other closed ${\rm SO}(r+1)$-orbits in the decomposition (\ref{eq:harmdec}).

\begin{prop0}\label{prop:ql} Let $q$ be a nondegenerate quadratic form, let $\ell$ be a linear form such that $\ell=0$ is  tangent to $Z(q)$.
Then
$${\rm Hess}(q^k\ell^h)=\left\{\begin{matrix}cq^{(r+1)(k-1)}\ell^{(r+1)h}&\textrm{\ if\ }k\ge 1\textrm{\ for\ }c=-2^{r-1}k^r(k+h)(2k+h-1),\\
0\textrm{\ if\ }k=0.\end{matrix}\right.$$
\end{prop0}
\begin{proof}
Let $P_\ell{\subseteq SO(V)}$ be the isotropy group acting on $V^\vee$ which fixes the hyperplane $Z(\ell)$ with equation $\ell=0$, as in Proposition \ref{prop:isotropy}.
The group $P_\ell$ acts on {$\PP(V)$} which is a prehomogeneous space with a dense orbit with complement the reducible divisor $Z(q)\cup Z(\ell)$.
The group  $P_\ell$  acts on each space $\PP(\sym^d (V^\vee))$.
{ Since $g\in \PP(\sym^d (V^\vee))$ is a fixed point if and only if $Z(g)\subset \PP(V)$ is an invariant subset, it follows from Proposition \ref{prop:isotropy} that} the only fixed points are $q^k\ell^h$ for $2k+h=d$.
Since the Hessian is a $P_\ell$-equivariant map, 
  ${\rm Hess}(q^k\ell^h)$ must be $cq^{s}\ell^{t}$ for some integer $s$, $t$ such that $2s+t=(r+1)(2k+h-2)$ and some scalar $c$.
We may set $q=x_0x_1+x_2^2+\cdots + x_{r}^2$, $\ell =x_0$. 

 We can now see which monomials of the type  $x_0^{\alpha}x_1^{\beta}$ occur in the hessian of $q^k\ell^h$. 
Let's take notice of these monomials occurring in each entry of the Hessian matrix, as follows:
{\footnotesize
$$\begin{pmatrix}(h+k)(h+k-1)x_0^{h+k-2}x_1^k&k(k+h)x_0^{h+k-1}x_1^{k-1}&0&\ldots&0\\
k(k+h)x_0^{h+k-1}x_1^{k-1}&k(k-1)x_0^{h+k}x_1^{k-2}&0&\ldots&0\\
0&0&2kx_0^{h+k-1}x_1^{k-1}&0&0\\
0&0&0&\ddots&0\\
0&0&0&0&2kx_0^{h+k-1}x_1^{k-1}\end{pmatrix}$$
}
Computing the determinant, we see that the only monomial containing $x_0, x_1$ is  $x_0^{(h+k-1)(r+1)}x^{(k-1)(r+1)}$,
with nonzero scalar coefficient if $k\ge 1$. The computation of $c$ is straightforward. This concludes the proof. \qedd \end{proof}\medskip

\begin{prop0}\label{prop:taylor_qk}
Let $\ell$ be a $q$-isotropic linear form, i.e., the hyperplane with equation $\ell=0$ is tangent to the quadric $Z(q)$.  We have 
$$\mathrm{Hess}\left(q^k+\epsilon q^{k-m}\ell^{2m}\right) = 
c_0q^{(r+1)(k-1)}+\epsilon c_1q^{(r+1)(k-1)-m}   \ell^{2m}+\cdots$$ 
where
$$\left\{\begin{array}{ll} 
c_0=2^{r-1}k^{{r+1}}(1-2k)\\
c_1=2^{r}k^{r}(2k-1)\left(2m^2+m(r-1)-k(r+1)\right).\end{array}\right.$$
\end{prop0}
\begin{proof}
As in the proof of Proposition \ref{prop:ql}, the terms of Taylor expansion are $P_\ell$- invariants, hence they are linear combinations of terms  of the form  $q^h\ell^j$
for integers $h, j$. The $0$-term has been computed in  Proposition \ref{prop:ql}.
Exactly as in the proof of Proposition \ref{prop:ql} we take note of the monomials containing  only $x_0, x_1$, in each entry of the Hessian matrix. 
The monomials containing  only $x_0$, $x_1$  fill up the following blocks of the Hessian matrix
$$\begin{pmatrix}A&0\\0&B\end{pmatrix}$$
where $A$ is the following $2\times 2$ matrix
{\footnotesize
$$\begin{pmatrix}k(k-1)x_0^{k-2}x_1^k+\epsilon(k+m)(k+m-1)x_0^{k+m-2}x_1^{k-m}&k^2x_0^{k-1}x_1^{k-1}+\epsilon(k^2-m^2)x_0^{k+m-1}x_1^{k-m-1}\\
k^2x_0^{k-1}x_1^{k-1}+\epsilon(k^2-m^2)x_0^{k+m-1}x_1^{k-m-1}&k(k-1)x_0^{k}x_1^{k-2}+\epsilon(k-m)(k-m-1)x_0^{k+m}x_1^{k-m-2}\end{pmatrix}$$
}
and $B$ is the identity matrix of size $r-1$ multiplied by 
$$2kx_0^{k-1}x_1^{k-1}+\epsilon 2(k-m)x_0^{k+m-1}x_1^{k-m-1}$$

By expanding the determinant, according to Proposition \ref{prop:ql} there is a term
$$c_0x_0^{(r+1)(k-1)}x_1^{(r+1)(k-1)}$$
corresponding to $c_0q^{(r+1)(k-1)}$. The $\epsilon$-term is 
$$c_1\epsilon x_0^{(r+1)(k-1)+m}x_1^{(r+1)(k-1)-m}$$ corresponding to
$\epsilon c_1q^{(r+1)(k-1)-m}   \ell^{2m}$.
\end{proof}
\vskip 0.4cm

We have the following consequence
\begin{prop0}\label{prop:maxrank} Let $\dim (V)=r+1$, $d=2k$.  Assume that for every $m$ such that $1\le m\le k$ we have
$2m^2+m(r-1)-k(r+1)\neq 0$. Then the Hessian map has maximal rank at $q^k$.
\end{prop0}
\begin{proof} The  differential  of the Hessian map at $q^k$  is a linear map between two representation spaces for ${\rm SO}(r+1)$, that is invariant for the Lie algebra of ${\rm SO}(r+1)$ . 
 By (\ref{eq:harmdec}) the  domain is $\oplus_{i=1}^{k} q^{k-i}\H_{2i}$.
Note that  the $i=0$ summand does not appear since the point corresponding to $q^{k}$ is a generator of that summand when seen in $\sym^d(V) $,
but here we are considering the tangent space at the corresponding point in the projective space   $\PP(\sym^d(V)) $. 
In the same way the target space is $\oplus_{i=1}^{(r+1)(k-1)} q^{(r+1)(k-1)-i}\H_{2i}$. Note that each summand in the  domain
has a corresponding isomorphic space in the target corresponding to the same index.

 So the differential of the Hessian map at $q^k$
sends a form  $ q^{k-i}\ell^{2i}\in q^{k-i}\H_{2i}$ to  a  non--zero scalar multiple of $q^{(r+1)(k-1)-i}   \ell^{2i}\in q^{(r+1)(k-1)-i}\H_{2i}$
{ by Proposition \ref{prop:taylor_qk}}.
Since this linear map is ${\rm SO}(r+1)$-equivariant,  
extending it by linearity it coincides, by  Schur's Lemma, with a scalar multiple of the identity (up to isomorphism) in each irreducible harmonic subspace.
 At level of polynomials, the differential of $h_{d,r}$ at $q^k$, restricted to the summand  $q^{k-i}\H_{2i}$, is the multiplication by a  non--zero scalar multiple of $q^{r(k-1)-1}$.
The scalar of each summand depend on $i$, so that  the  differential of $h_{d,r}$ at $q^k$  has a diagonal block structure padded by zeros.
\qedd \end{proof}

\begin{rem0}\label{rem:bundles} {\rm
As in the proof of Proposition \ref{prop:ql}, we set $q=x_0x_1+x_2^2+\cdots + x_{r}^2$, $\ell =x_0$, a $q$-isotropic linear form.
Every matrix in $P_\ell$ fixes $\ell$, hence it has $(1,0,\ldots, 0)^t$ as eigenvector. Then every $A\in P_\ell$ has its entries satisfying $a_{i1}=0$ for $i\ge 2$ and $a_{11}\neq 0$. 
The one dimensional representation of $P_\ell$ given by $A\mapsto a_{11}^{-1}$, where $a_{11}$ is the $(1,1)$-entry (necessarily nonzero) of $A$, is called $L$.
We see now that this representation corresponds to a basic line bundle.

There is a well known equivalence of categories  between $P_\ell$-modules and ${\rm SO}(r+1)$-equivariant vector bundles on the smooth $(r-1)$-dimensional quadric
 $\mathbb Q_{r-1}={\rm SO}(r+1)/P_\ell$ (see \cite{Kap, BK}). In this latter category the morphisms are ${\rm SO}(r+1)$-equivariant morphisms.
In a nutshell, given $\rho\colon P_\ell\to {\rm GL}(r)$, it is defined an action of $P_\ell$ on  ${\rm SO}(r+1)\times \CC^r$  given by
$p\cdot (g,v) =(gp,\rho(p^{-1})v)$ for $p\in P_\ell$,  $g\in {\rm SO}(r+1)$,  $v\in\CC^r$, and the orbit space has a natural vector bundle structure
 $E_{\rho}$ with fibers $\CC^r$  on the variety 
 $\mathbb Q_{r-1}$. In other words, $E_{\rho}$ is obtained by the principal bundle ${\rm SO}(r+1)\to {\rm SO}(r+1)/P_\ell=\mathbb Q_{r-1}$ via the homomorphism $\rho$. 

Under this equivalence of categories, the representation $L$ corresponds to the line bundle $\O(1)$ on  $\mathbb Q_{r-1}$, 
 and more generally
the representation $L^{\otimes t}$ corresponds to $\O(t)$. This can be seen since the space of sections  $H^0(\mathbb Q_{r-1}, E_{\rho})$ is identified
with $\left\{f\colon {\rm SO}(r+1)\to\CC^r| f(gp)=\rho(p^{-1})f(g)\right\}$ with the natural ${\rm SO}(r+1)$-action given by $(g\cdot f)(g_1):=f(g^{-1}g_1)$. 
Then a basis of $H^0(\mathbb Q_{r-1}, L)$ is given by $f_i(A)=a_{i1}$ for $i=0,\ldots, r$, which identifies $H^0(\mathbb Q_{r-1},L)\cong V$ as ${\rm SO}(r+1)$-modules, hence $L=\O(1)$ as  wanted.
}
\end{rem0}

Let $\ell$ be a  non--zero $q$-isotropic linear form, and let $d=2k+1$.  The differential of the Hessian map $h_{d,r}$ at $q^k\ell$ is a linear map invariant for the isotropy group $P_\ell\subset {\rm SO}(r+1)$,
consisting of $g\in {\rm SO}(r+1)$ such that $g(\ell)=\ell$. The group  $P_\ell$ is well known  to be  parabolic \cite[\S 23.3]{FH} but  not reductive for $r\ge 2$.

The differential of $h_{d,r}$ at $q^k\ell$ is as follows 

\begin{equation} \label{eq:differential0} h'_{d,r}\colon\left(T(\PP(\sym^{d}V))\right)_{q^k\ell}\longrightarrow \left(T(\PP(\sym^{(r+1)(d-2)}V))\right)_{q^{(r+1)(k-1)}\ell^{r+1}}.
\end{equation} 
where we used Proposition \ref{prop:ql}. We study now the rank of $h'_{d,r}$. We will show in  Proposition \ref{prop:maxrankql} that $h'_{2k+1,r}$ has maximal rank for $k\ge 2$, 
with  some  numerical assumption on $k$, $r$,
exactly as we did in the even case in
Proposition \ref{prop:maxrank}. In order to linearize this map, it is convenient to consider the induced map
$$\left(T(\PP(\sym^{d}V))\right)_{q^k\ell}\longrightarrow \left(h_{d,r}^*(T(\PP(\sym^{(r+1)(d-2)}V)))\right)_{q^k\ell}$$ where now both  spaces
are based at the same point $q^k\ell$.

The rank of the differential $h'_{d,r}$ in \eqref{eq:differential0} is  one less than  the rank of the corresponding linear map between vector spaces, namely the map

\begin{equation} \label{eq:differentialv} \bigoplus_{j=0}^k \H_{d-2j}
\longrightarrow \bigoplus_{i=0}^{ (r+1)(k-1)} \H_{(r+1)(d-2)-2i}\otimes L^r\end{equation}
where the exponent $r$ of $L$  is  computed by the fact that $h_{d,r}(q^k\ell ) = q^{(r+1)(k-1)}\ell^{r+1}$,
so that $q^k\ell \in \H_1$ goes to $ q^{(r+1)(k-1)}\ell^{r+1}\in \H_{r+1}\otimes L^r$,  and we see that the factor $\ell$  is multiplied by $\ell^r\in H^0(\O(r))$
and the line bundle $\O(r)$ corresponds to $L^r$ (see Remark \ref{rem:bundles}). Note that the map $h'_{d,r}$ in (\ref{eq:differential}) is injective if and only if the map in
(\ref{eq:differentialv}) is injective.

Comparing with the proof of  Proposition \ref{prop:maxrank}, we note that since $P_\ell$ is not reductive,
the summands of  (\ref{eq:differentialv}) are no  longer  irreducible for the action of $P_\ell$
and we can no longer use Schur's Lemma. 
For example $\H_1=V$ has the { nonsplitting} filtration
$0\subset\langle\ell\rangle\subset\langle\ell\rangle^\perp \subset V$
where  the consecutive quotients are irreducible, and  similar considerations hold  for all $\H_r$. { The source and the target of  (\ref{eq:differentialv}) both split for the action of the reductive part of $P_\ell$ in summands having multiplicities which are difficult to control, so that we find more efficient to consider the map (\ref{eq:differentialv})
as a map of $P_\ell$-modules.}

 A description of $P_\ell$-modules in terms of quiver representations is exposed in \cite{OR}, but here we proceed in a more elementary way. Let us define first the following linear maps.
\begin{defn0}\label{def:pijr}  Let  $i,j, k$ be nonnegative integers, such that $|j-i|\le k\le i+j$, $i+k-j$ being even.
The linear  map $P_{i,j}^k\colon\H_i\to\H_j\otimes L^k$ is defined by
$$P_{ij}^k(h):=(h\ell^k)_j$$ where $(h\ell^k)_j$ is the $j$-th harmonic summand of $h\ell^k$, according to Propoposition \ref{prop:harmdec}. 
 Notice that $j=k+i\ge 0$ is the maximal value of $j$ for which  $P_{ij}^k$ is nonzero, for  given $i$, $k$, { since the degree of $h\ell^k$ is $k+i$.}
\end{defn0}
\begin{prop0}\label{prop:pij_inj}
The map $P_{ij}^k$ is $P_\ell$-equivariant and it is injective for $j=k+i\ge 0$  (in this case the first inequality in Definition \ref{def:pijr} becomes an equality). 
\end{prop0}
\begin{proof}
The equivariance is clear since $\ell$ is a $P_\ell$-invariant  function (up to scalar multiples)
and the projection on the harmonic summands are even ${\rm SO}(V)$-invariant,  see Remark \ref{rem:bw}. Let $j=k+i$ and assume $g\in {\rm Ker}( P_{ij}^{j-i})$. 
Then $g\ell^{j-i}$ has degree $j$
and has vanishing harmonic part, hence it is divided by $q$. It follows that $g$ is divided by $q$ which implies $g=0$ since $g$ is harmonic.\qed
\end{proof}

The following Lemma is a well known result from Representation Theory.
\begin{prop0}\label{prop:tensorharm}
Let $i\le j$. As ${\rm SO}(r+1)$-modules, we have the following decompositions
\begin{itemize}
\item{}
If $r\ge 3$ then $\H_i\otimes \H_j = \left(\oplus_{p=0}^i\H_{i+j-2p}\right)\oplus T$ where $T$ is the direct sum of certain irreducible  ${\rm SO}(r+1)$-modules,
not isomorphic to any $\H_q$.

\item{} If $r=2$ then $\H_i\otimes\H_j=\oplus_{k=j-i}^{j+i}\H_k$ (note here there is no parity condition among $i$, $j$ and $k$).

\end{itemize}
\end{prop0}
\begin{proof}  It follows from the Littlewood-Richardson rule for the orthogonal group,
which is exposed in full generality in \cite{Lit}. The particular cases treated in \cite{KT} are enough for our purposes, see  \cite [Example (2), p. 510]{KT}.  Actually this example
is exposed for the symplectic group, but anything  applies also  in the case of the orthogonal group, see \cite[p. 509]{KT}. 
In the case $r=2$, the module  $\H_d$  has dimension $2d+1$ and corresponds to $H^0(\PP^1,\O(2d))$, in this case the result is classical and it is sometimes attributed to Clebsch-Gordan
(see \cite[exerc. 11.11]{FH} for a modern reference).\qed
\end{proof}\smallskip

Although we do not use it, it seems worth to state the following interesting consequence of Proposition \ref{prop:tensorharm}.
\begin{coro0}  Let $r\ge 2$. Let $f_i\in\H_i$, $f_j\in\H_j$. Then the $k$-th harmonic summands of the product $f_if_j$ are nonzero
only for $|j-i|\le k\le i+j$, $i+j-k$ even.
\end{coro0}
\begin{proof} The map which sends the pair $(f_i, f_j)$ the $k$-th harmonic summand of the product $f_if_j$
is a ${\rm SO}(r+1)$-equivariant map $\H_i\otimes\H_j\to\H_k$. The parity condition follows from the harmonic decomposition and it is necessary also for $r=2$. 
In other words, among the summands of the tensor product  $\H_i\otimes\H_j$ which are listed  in the case $r=2$ of Proposition \ref{prop:tensorharm},
only the ones satisfying the condition $i+j-k$ even appear
as a $k$-th harmonic summand of $f_if_j$.
\end{proof}

\begin{prop0}\label{prop:LitRich}
Let $r\ge 2$. Every nonzero $P_\ell$-equivariant map from $\H_i$ to $\H_j\otimes L^k$ coincides with a scalar multiple of $P_{i,j}^k$ as in Definition \ref{def:pijr}.
In particular one has $|j-i|\le k\le i+j$, $i+j-k$ is even.
\end{prop0}
\begin{proof} We come back to the equivalence of categories  between $P_\ell$-modules and ${\rm SO}(r+1)$-equivariant vector bundles on  $\mathbb Q_{r-1}={\rm SO}(r+1)/P_\ell$ 
sketched in Remark \ref{rem:bundles}.
The injective map $P_{ij}^{j-i}$ of Proposition \ref{prop:pij_inj} corresponds to the injective  map
$$H^0(\mathbb Q_{r-1},\O(i))\otimes \O_{\mathbb Q_{r-1}}\to H^0(\mathbb Q_{r-1},\O(j))\otimes\O(j-i)$$
This map is natural:  indeed it corresponds to the multiplication of sections of $\O(i)$ with sections of $\O(j-i)$, which gives sections of $\O(j)$. 
This map is ${\rm SO}(r+1)$-equivariant  and corresponds, taking the $H^0$  functor, to the injection
$\H_i\to \H_j\otimes \H_{j-i}$.
Consider now any nonzero $P_\ell$-equivariant map from $\H_i$ to $\H_j\otimes L^k$. By taking the $H^0$  functor, such a map corresponds to a ${\rm SO}(r+1)$-equivariant map $\H_i\to \H_j\otimes \H_{k}$ (see lemma \ref {lem:isoquad}) and it exist exactly for the values of $k$ considered in Definition \ref{def:pijr}
and Proposition \ref{prop:tensorharm}.
\qed
\end{proof}

\begin{rem0} {\rm In the case $r=2$ we have that the Lie algebra of ${\rm SO}(3)$  is ${\mathfrak sl(2)}$ and the tensor product contains more summands,
look at Proposition \ref{prop:tensorharm}, where for $r=2$ the parity condition disappear. Apparently we get more maps  than the ones stated in  Proposition \ref{prop:LitRich},
but the statement there is actually correct.
There are indeed maps of homogeneous bundles on ${\mathbb Q}_1\cong \PP^1$  that do not come from representations of $P_\ell$.
 In this case ${\mathbb Q}_1\cong \PP^1$  is embedded in the plane $\PP(\H_1)$
by twice the generator of ${\rm Pic}(\PP^1)$. Consider the  universal covering ${\rm SL(2)}\rig{\pi}{\rm SO}(3)$. The generator of ${\rm Pic}(\PP^1)$ is homogeneous as well,
but it comes from a representation of the parabolic group $\pi^{-1}(P_\ell)\subset {\rm SL}(2)$, so it does not come from a representation of $P_\ell$.

}
\end{rem0}

\begin{prop0}\label{prop:taylor_qkl}
Let $\ell$ be a $q$-isotropic linear form.  We have 
 $$\mathrm{Hess}\left(q^k\ell+\epsilon q^{k-m}\ell^{2m+1}\right) = 
c_0q^{(r+1)(k-1)}\ell^{r+1}+\epsilon c_1q^{(r+1)(k-1)-m}   \ell^{2m+r+1}+\cdots$$ 
where
$$\left\{\begin{array}{ll} 
c_0=-2^{r}k^{r+1}(k+1)\\
c_1=2^{r}k^{r}\left(m^2(2k+1)+m(rk+r-k)-k(k+1)(r+1)\right).\end{array}\right.$$
\end{prop0}
\begin{proof}
As in the proof of Proposition \ref{prop:ql}, the terms of Taylor expansion are $P_\ell$- invariants, hence they are linear combinations of terms $q^h\ell^j$
for integers $h, j$. The $0$-term has been computed in  Proposition \ref{prop:ql}.
Exactly as in the proof of Proposition \ref{prop:ql} we take note of the monomials containing just $x_0, x_1$, in each entry of the Hessian matrix,
that  fill  up  the following blocks 
$$\begin{pmatrix}A&0\\0&B\end{pmatrix}$$
where $A$ is the following $2\times 2$ matrix
{\footnotesize
$$\begin{pmatrix}k(k+1)x_0^{k-1}x_1^k+\epsilon(k+m)(k+m+1)x_0^{k+m-1}x_1^{k-m}&k(k+1)x_0^{k}x_1^{k-1}+\epsilon(k-m)(k+m+1)x_0^{k+m}x_1^{k-m-1}\\
k(k+1)x_0^{k}x_1^{k-1}+\epsilon(k-m)(k+m+1)x_0^{k+m}x_1^{k-m-1}&k(k-1)x_0^{k+1}x_1^{k-2}+\epsilon(k-m)(k-m-1)x_0^{k+m+1}x_1^{k-m-2}\end{pmatrix}$$
}
and $B$ is the identity matrix of size $(r-1)$ multiplied by 
$$2kx_0^{k}x_1^{k-1}+\epsilon 2(k-m)x_0^{k+m}x_1^{k-m-1}$$

By expanding the determinant, according to Proposition \ref{prop:ql}, there is a term
$$c_0x_0^{(r+1)k}x_1^{(r+1)(k-1)}$$
corresponding to $c_0q^{(r+1)(k-1)}\ell^{r+1}$. The $\epsilon$-term is 

$$c_1\epsilon x_0^{(r+1)k+m}x_1^{(r+1)(k-1)-m}$$ corresponding to
$\epsilon c_1q^{(r+1)(k-1)-m}   \ell^{2m+r+1}$.\qed
\end{proof}

\begin{prop0}\label{prop:taylor_qkl2}
Let $\ell$ be a $q$-isotropic linear form.  We have the formula
  $$\mathrm{Hess}\left(q^{k-1}\ell^2+\epsilon q^{k-m}\ell^{2m}\right) = 
c_0q^{(r+1)(k-2)}\ell^{2(r+1)}+\epsilon c_1q^{(r+1)(k-2)+1-m}   \ell^{2m+2r}+\cdots$$ 
where
$$\left\{\begin{array}{ll} 
c_0=-2^{r-1}(k-1)^{r}(k+1)(2k-1)\\
c_1=2^{r-1}(k-1)^{r-1}(2k-1)\left(2km^2+m(rk+r-5k+1)-k(k(r+1)+r-3)\right).\end{array}\right.$$
\end{prop0}
\begin{proof}
As in the proof of Proposition \ref{prop:ql}, the terms of Taylor expansion are $P_\ell$- invariants, hence they are linear combinations of terms $q^h\ell^j$
for integers $h, j$. The $0$-term has been computed in  Proposition \ref{prop:ql}.
Exactly as in the proof of Proposition \ref{prop:ql} we take note of the monomials containing just $x_0, x_1$, in each entry of the Hessian matrix, that fill up the following blocks
$$\begin{pmatrix}A&0\\0&B\end{pmatrix}$$
where $A$ is the  $2\times 2$ matrix
{\footnotesize
$$\begin{pmatrix}k(k+1)x_0^{k-1}x_1^{k-1}+\epsilon(k+m)(k+m-1)x_0^{k+m-2}x_1^{k-m}&(k^2-1)x_0^{k}x_1^{k-2}+\epsilon(k-m)(k+m)x_0^{k+m-1}x_1^{k-m-1}\\
(k^2-1)x_0^{k}x_1^{k-2}+\epsilon(k-m)(k+m)x_0^{k+m-1}x_1^{k-m-1}&(k-1)(k-2)x_0^{k+1}x_1^{k-3}+\epsilon(k-m)(k-m-1)
x_0^{k+m}x_1^{k-m-2}\end{pmatrix}$$
}
and $B$ is the identity matrix of size $(r-1)$ multiplied by 
$$2(k-1)x_0^{k}x_1^{k-2}+\epsilon 2(k-m)x_0^{k+m-1}x_1^{k-m-1}$$

By expanding the determinant, according to Proposition \ref{prop:ql} there is a term

$$c_0x_0^{(r+1)k}x_1^{(r+1)(k-1)}$$
corresponding to $q^{(r+1)(k-1)}\ell^{r+1}$. The $\epsilon$-term is 

$$\epsilon c_1x_0^{(r+1)k+m}x_1^{(r+1)(k-1)-m}$$ corresponding to
$\epsilon c_1q^{(r+1)(k-1)-m}   \ell^{2m+r+1}$. \qed
\end{proof}

\begin{prop0}\label{prop:maxrankql} Let $d=2k+1$, $r\ge 2$.  Assume that for every $m$ such that 
$0\le m \le k$ we have
 $m^2(2k+1)+m(rk+r-k)-k(k+1)(r+1)\neq 0$. Then 
 the differential of $h_{d,r}$  at $q^k\ell$ is injective.
\end{prop0}
\begin{proof}  The differential of $h_{d,r}$ at $q^k\ell$ as in (\ref{eq:differentialv}) takes the summand  $\mathcal H_{d-2j}$ in the domain to several summands $\H_p$ in the target,
for $p\le d-2j+r$, as it follows from Proposition  \ref{prop:LitRich}.  

 By using Proposition \ref{prop:LitRich}, the map induced from $\H_{d-2j}$ to the extreme summand $\H_{d-2j+r}$ is a nonzero scalar multiple of $P_{d-2j,d-2j+r}^r$ defined in Definition \ref{def:pijr}, hence it is injective.
 Indeed,  the fact that this map is nonzero follows for $j<k$ by computing
the derivative at $q^k\ell$ in the direction of  $q^{k-m}\ell^{2m+1}$ with the {\it same} linear form
$\ell$,
which makes the computation easier and it has been already performed in Proposition \ref{prop:taylor_qkl}. Here we use the numerical assumption in the hypothesis.
In particular, the differential of $h_{d,r}$ has a triangular block shape, with the diagonal blocks of maximal rank, hence it has maximal rank and it is indeed injective.

 In order to
prove that the restriction of $h_{d,r}$ is nonzero also on the summand $\H_1/\langle q^{k}\ell\rangle$, for $j=k$,  we fix $q=x_0x_1+\sum_{i\ge 2}x_i^2$ and we consider the Hessian of $q^k x_0+\epsilon q^{k}x_2$ which is 
$$-2^rk^{r+1}(k+1)q^{(r+1)(k-1)}\left(x_0^{r+1}+\epsilon (r+1) x_0^{r}x_2+\cdots\right).$$ 
This computation can be performed, as in the proof of Proposition \ref{prop:ql}, in the case where $x_i=0$ for $i\ge 3$, so looking only at the monomials containing $x_0$, $x_1$, $x_2$.   In this computation $x_2$ may be replaced with any other linear form
$q$-orthogonal with $x_0$. The assertion  follows. \qed
\end{proof}

\begin{prop0}\label{prop:maxrankql2} Let $d=2k$, $r\ge 2$. Assume that for every $m$ such that 
$0\le m\le k$ we have $2km^2+m(rk+r-5k+1)-k(k(r+1)+r-3)\neq 0$. Then the differential of $h_{d,r}$ at  $q^{k-1}\ell^2$  is injective.
\end{prop0}
\begin{proof} The argument is analogous to the one of  Proposition \ref{prop:maxrankql} and we just sketch it.
The differential of $h_{d,r}$ at $q^{k-1}\ell^2$ is the  $P_\ell$-equivariant map

\begin{equation} \label{eq:differential} h'_{d,r}\colon\left(T(\PP(\sym^{d}V))\right)_{q^{k-1}\ell^2}\longrightarrow \left(T(\PP(\sym^{(r+1)(d-2)}V))\right)_{q^{(r+1)(k-2)}\ell^{2(r+1)}}.
\end{equation} 
which is injective if and only if the corresponding linear map

\begin{equation} \label{eq:differentialv2} \bigoplus_{j=0}^k \H_{d-2j}
\longrightarrow \bigoplus_{i=0}^{ (r+1)(k-1)} \H_{(r+1)(d-2)-2i}\otimes L^{2r}\end{equation}
is injective.

We continue exactly as in Proposition \ref{prop:maxrankql}, by using Propositions \ref{prop:LitRich} and \ref{prop:taylor_qkl2}.
 In order to
prove that the restriction  of $h_{d,r}$ is nonzero also on the summand $\H_2/\langle q^{k-1}\ell^2\rangle$,   we fix $q=x_0x_1+\sum_{i\ge 2}x_i^2$ and  we consider the Hessian of $q^{k-1} x_0^2+\epsilon q^{k-1}x_0x_2$ which is 
$$-2^{r-1}(k-1)^{r}(k+1)(2k-1)q^{(r+1)(k-2)}\left(x_0^{2(r+1)}+\epsilon (r+1) x_0^{2r+1}x_2+\cdots\right).$$ 
\qed \end{proof}

\begin{rem0} {\rm 
The parabolic subgroup $P_\ell$ used in the proof of Proposition \ref{prop:ql} appears implicitly in
\cite[\S 1, Example 4]{Dol}, where it is considered the homaloidal polynomial $x_0(x_0x_2+x_1^2)$ as relative invariant for $P_\ell\subset {\rm SO}(3)$.
In the same way, the  Perazzo cubic $3$-fold (see for example \cite[Theorem 7.6.7 (iv)]{Ru}), is a relative invariant for the  parabolic subgroup
$P\subset {\rm SL}(3)$ of linear transformations fixing a point, which is $6$-dimensional. The Perazzo cubic $3$-fold is isomorphic to the zero locus of
the symmetric determinant 
$$\begin{vmatrix}0&x_0&x_1\\x_0&x_2&x_3\\x_1&x_3&x_4\end{vmatrix}$$
 and  this determinantal expression shows it can be interpreted as the variety of singular conics
 passing through a point $p\in\PP^2$ . 
This remark completes the description of the isotropy group of the Perazzo cubic $3$-fold begun in \cite[Example 4.3]{CNOS}.}
\end{rem0}

\begin{rem0}\label{rem:contracts} {\rm Proposition \ref{prop:ql} gives an obstruction to extend the present technique to the cubic case in higher dimension.
Indeed, since  ${\rm Hess}(q\ell)=-2^{r+1}\ell^{r+1}$, the Hessian map $h_{3,r}$ contracts the variety $E_\ell=\{q\ell | q\textrm{\ is tangent to\ }\ell\}$
(which is a degree $r$ hypersurface)
to a point and it has not maximal rank at $q\ell$.}
\end{rem0}

We prove now our main criterion of birationality by using the group action.
\begin{thm0}\label{thm:criterion}
Let $G_{d,r}\subset\PP(\sym^d (V))\times\PP(\sym^{(d-2)(r+1)}(V))$ be the closure of the graph
$\left\{\left(f,h_{d,r}(f)\right), \textrm{ for any }f\textrm{ such that  }h_{d,r}(f)
\textrm{ is defined}\right\}$.
Assume there exists a nondegenerate quadratic form $q$ and a linear form $\ell$  such that the hyperplane $\ell=0$ is
tangent to the quadric $Z(q)$ and  such that  one of the following three assumptions hold 
\begingroup
\allowdisplaybreaks
\begin{align}\label{eq:notin}
\begin{split}
(\ell^d,q^{(r+1)(k-1)})\notin G_{2k,r}\textrm{\ if\ }d=2k\textrm{\ is even} \\ 
\textrm{and moreover assume that for every\ }m\textrm{ such that\ } 1\le m\le k\\
\textrm{ we have\ }
2m^2+m(r-1)-k(r+1)\neq 0;\end{split}\end{align}
\endgroup

\begingroup
\allowdisplaybreaks
\begin{align}\label{eq:notinodd}
\begin{split}
(\ell^d,q^{(r+1)(k-1)}\ell^{r+1})\notin G_{2k+1,r}\textrm{\ if\ }d=2k+1\textrm{\ is odd} \\ 
\textrm{and moreover assume that for every\ }m\textrm{ such that\ } 0\le m\le k\\
\textrm{ we have\ }
m^2(2k+1)+m(rk+r-k)-k(k+1)(r+1)\neq 0;\end{split}\end{align}
\endgroup

\begingroup
\allowdisplaybreaks
\begin{align}\label{eq:notineven2}
\begin{split}
(\ell^d,q^{(r+1)(k-2)}\ell^{2(r+1)})\notin G_{2k,r}\textrm{\ if\ }d=2k\textrm{\ is even} \\ 
\textrm{and moreover assume that for every\ }m\textrm{ such that\ } 0\le m\le k\\
\textrm{ we have\ }
2km^2+m(rk+r-5k+1)-k(k(r+1)+r-3)\neq 0.\end{split}\end{align}
\endgroup

Then $h_{d,r}$ is birational onto its image.
\end{thm0}
\begin{proof}
The group $G={\rm SO}(V)$ leaves $G_{d,r}$ invariant and the projection
 on the second factor $\pi\colon G_{d,r}\to\PP(\sym^{(d-2)(r+1)}(V))$ is $G$-equivariant.
We consider first the case $d$ even.
We have to prove that $\pi$ is generically injective and since it is a morphism between projective varieties
it is enough to prove that the scheme-theoretic fiber  of some point is a unique point.  Under the assumption (\ref{eq:notin}), we may want to prove that   
$\pi^{-1}(q^{(r+1)(k-1)})$ consists of the single pair
$(q^k,q^{(r+1)(k-1)})$ (compare with Proposition \ref{prop:hescon}).
The fiber $\pi^{-1}(q^{(r+1)(k-1)})$ is $G$-invariant. Since $(q^k,q^{(r+1)(k-1)})$ is a connected component
of the fiber by Proposition \ref{prop:maxrank}, if the fiber is larger it must contain another closed $G$-orbit.
These closed orbits have the form 
$$\left\{(q^a\ell ^b,q^{(r+1)(k-1)})\textrm{ for }\ell \textrm{ isotropic}\right\}$$
 for $a$ a non--negative integer and $b$ a positive integer, such that $2a+b=2k$. 
 The cases $a>0$ do not belong to $G_{d,r}$ by Proposition
\ref{prop:ql}. The case $a=0, b=2k$ does not belong to $G_{d,r}$ by the assumption (\ref{eq:notin}).

Similarly,  under the assumption (\ref{eq:notineven2}), given any isotropic linear form $\ell$, we may want to prove that   the fiber 
$\pi^{-1}(q^{(r+1)(k-2)}\ell^{2(r+1)})$ consists of the single pair $(q^{k-1}\ell^2,\\ q^{(r+1)(k-2)}\ell^{2(r+1)})$ (compare with Proposition \ref{prop:ql}).
Consider the algebraic set $\Sigma$ of all forms of the type $q^{(r+1)(k-2)}\ell^{2(r+1)}$ with $\ell$ isotropic. The set $\Sigma$ is $G$--invariant, and therefore also its counterimage $\Sigma'$ via $\pi$ is $G$--invariant. By Proposition \ref {prop:maxrankql2}, there is a connected component $\Sigma''$ of $\Sigma'$, that consists of all pairs of the form $(q^{k-1}\ell^2,\\q^{(r+1)(k-1)}\ell^{2(r+1)})$, with $\ell$ isotropic. We want to prove that $\Sigma'=\Sigma''$. We argue by contradiction. If this is not the case, then there is some closed orbit contained in $\Sigma'\setminus \Sigma''$. Such a closed orbit, different from $\Sigma''$, could only be of the form
$$\left\{(q^a\ell ^b,q^{(r+1)(k-2)}\ell^{2(r+1)})\textrm{ for }\ell \textrm{ isotropic}\right\}$$
with $a,b$ non--negative integers such that $2a+b=2k$ and $a\neq k-1$. 
 The cases $a>0$, $a\neq k-1$ give a contradiction by Theorems \ref {prop:hescon} and \ref {prop:ql}. The case $a=0$ and $b=2k$ is excluded by \eqref {eq:notineven2}. 
 This concludes the proof of the case $d$ even. 
 
 The case $d$ odd,  under the assumption (\ref{eq:notinodd}),  is analogous by using Proposition \ref{prop:maxrankql}.\qedd
\end{proof}

\begin{rem0} {\rm A linear form $\ell$  defines a hyperplane $\ell=0$ that  is tangent to the nondegenerate quadric $Z(q)$ if and only if $\ell$ is an isotropic point for the dual quadric, and then we have called $\ell$ a $q$-isotropic linear form. Since $G_{d,r}$ is ${\rm SL}(V)$-invariant, if $(\ell^d,q^{(r+1)(k-1)})\in G_{d,r}$ for a nondegenerate quadratic form $q$
and a $q$-isotropic linear form $\ell$ (see  (\ref{eq:notin}) ), then $(\ell^d,q^{(r+1)(k-1)})\in G_{d,r}$ for any nondegenerate quadratic form $q$
and any $q$-isotropic linear form $\ell$.}
\end{rem0}
In the next section we will investigate the assumptions of Theorem \ref{thm:criterion}. We will prove in Corollary \ref{coro:indet}
that (\ref{eq:notin}) is satisfied 
when $r=2$ and $k\ge 2$  ($d$ even) and that  (\ref{eq:notinodd}) is satisfied 
when $r=2$ and $k\ge 3$  ($d$ odd) . Note that for $r=1$ this strategy does not work,
in \cite[Prop. 2.5]{CO} we proved that for $r=1$ and for any linear forms $x$, $y$ then $(x^d,(xy)^{d-2})\in G_{d,1}$,
so that we used a different approach to prove the birationality of $h_{d,1}$.
It is unclear if the assumptions (\ref{eq:notin}) and (\ref{eq:notinodd})  hold for $r\ge 3$.

\section{The indeterminacy values assumed by the Hessian map for ternary forms}
The goal of this section is to show that the assumptions of Theorem \ref{thm:criterion}, regarding the membership of certain elements to the closure of the graph of the Hessian map,  are satisfied for ternary forms, i.e.,  when $r=2$.
Recall that the Hessian map is not defined at forms which are cones but, approaching such forms, the Hessian map (as all the rational maps) can assume some ``indeterminacy values".
In \cite[Prop. 2.5]{CO} we proved that, for $r=1$,  approaching the power $x_0^d$, the Hessian is a form of degree $2d-4$ which is divisible by $x_0^{d-2}$,
so these are the indeterminacy values of the Hessian map in the case of binary forms.
The main result of this section is Theorem \ref{thm:divides}, which shows an analogous property for ternary forms.
\vskip 0.4cm

Let $G_{d,2}\subset \PP(\sym^d(\CC^{3}))\times\PP(\sym^{3(d-2)}(\CC^{3}))$
be the closure of the graph of the Hessian map $h_{d, 2}$ as in Theorem \ref{thm:criterion}.
 Assume   $(x_0^{d}, r)\in G_{d,2}$  for some polynomial $r$. The point $(x_0^{d},r)$
  may be approximated by an algebraic one dimensional branch $(f(t),g(t))$ such that $g(t)={\rm Hess}(f(t))$  for small $t\neq 0$
 and such that $f(0)=x_0^{d}$. Considering the projection $f(t)$ to the first component, we get a Puiseux series 
 
\begin{equation}\label{eq:fPuiseux}f(t)=x_0^d+\sum_{i=1}^{+\infty} t^{\alpha_i}f_i\end{equation}
with $\alpha_i\in\QQ$, $0<\alpha_1<\alpha_2<\cdots$, the denominators of $\alpha_i$ are bounded
and the series converges in the Euclidean topology for $|t|<\delta$.
We may assume that $h_{d,2}(f(t))={\rm Hess}(f(t))$ is well defined for $0<|t|<\delta$.

The Hessian matrix of $f(t)$ is
$$\begin{pmatrix}d(d-1)x_0^{d-2}+\sum t^{\alpha_i}f_{i,00}&\sum t^{\alpha_i}f_{i,01}&\sum t^{\alpha_i}f_{i,02}\\
\sum t^{\alpha_i}f_{i,01}&\sum t^{\alpha_i}f_{i,11}&\sum t^{\alpha_i}f_{i,12}&\\
\sum t^{\alpha_i}f_{i,02}&\sum t^{\alpha_i}f_{i,12}&\sum t^{\alpha_i}f_{i,22}\end{pmatrix},$$
 where, of course, the indices  refer to differentiation.

Computing the determinant, by linearity on each row, we have in a neighborhood of $t=0$ 
\begin{align}\label{eq:puiseux}{\rm Hess}(f(t))=\\
\nonumber d(d-1)x_0^{d-2}\left(\sum_{i, j}t^{\alpha_i+\alpha_j}H_{12}(f_i,f_j)\right)+\\
\nonumber \sum_{i, j, k}t^{\alpha_i+\alpha_j+\alpha_k}H(f_i, f_j, f_k)\end{align}
where
$$H_{12}(f,g)=H_{12}(g,f)=\frac 12\left(f_{11}g_{22}-2f_{12}g_{12}+f_{22}g_{11}\right)=\frac 12\left|\begin{matrix}f_{11}&f_{12}\\
g_{21}&g_{22}\end{matrix}\right|+\frac 12\left|\begin{matrix}g_{11}&g_{12}\\
f_{21}&f_{22}\end{matrix}\right|$$
(see for example \cite[\S 4.3]{O12}) and
{\tiny \begin{align}\nonumber 6H(f,g,h)=&\\ \nonumber\left|\begin{matrix}f_{00}&f_{01}&f_{02}\\
g_{10}&g_{11}&g_{12}\\
h_{20}&h_{21}&h_{22}\end{matrix}\right|+&
\left|\begin{matrix}f_{00}&f_{01}&f_{02}\\
h_{10}&h_{11}&h_{12}\\
g_{20}&g_{21}&g_{22}\end{matrix}\right|+\left|\begin{matrix}g_{00}&g_{01}&g_{02}\\
f_{10}&f_{11}&f_{12}\\
h_{20}&h_{21}&h_{22}\end{matrix}\right|+
\left|\begin{matrix}g_{00}&g_{01}&g_{02}\\
h_{10}&h_{11}&h_{12}\\
f_{20}&f_{21}&f_{22}\end{matrix}\right|+
\left|\begin{matrix}h_{00}&h_{01}&h_{02}\\
f_{10}&f_{11}&f_{12}\\
g_{20}&g_{21}&g_{22}\end{matrix}\right|+\left|\begin{matrix}h_{00}&h_{01}&h_{02}\\
g_{10}&g_{11}&g_{12}\\
f_{20}&f_{21}&f_{22}\end{matrix}\right|\end{align}}
so that $H(f,f,f)={\rm Hess}(f)$. We set $H_{12}(f)=H_{12}(f,f)=f_{11}f_{22}-f_{12}^2$.
\begin{lemma0}\label{lem:heszero}Let $r=2$.
Assume $f_{11}f_{22}-f_{12}^2$ is identically zero.
Then 
\begin{itemize}
\item{}$f=x_0^d+x_0^{d-1}l(x_1,x_2)+\sum_{i= 2}^dx_0^{d-i}c_im(x_1,x_2)^i$ with $l$, $m$ linear forms, $c_i$ scalars.
\item{}${\rm Hess}(f)$ is divisible by $x_0^{2d-4}$. Moreover ${\rm Hess}(f)$ vanishes if and only if  $l$ and $m$ are proportional. \end{itemize}
\end{lemma0}
\begin{proof} Let $f=\sum_{i=0}^df^{(i)}(x_1,x_2)x_0^{d-i}$ be the $x_0$-expansion of $f$,  with $f^{(i)}(x_1,x_2)$ a form of degree $i$, for $0\leq i\leq d$.  The first two summands $f^{(0)}x_0^d+f^{(1)}x_0^{d-1}$ do not contribute to $f_{11}f_{22}-f_{12}^2$, hence the assumption has an influence only on the summands for $2\leq i \leq d$.
This explains the different behaviour of the summands in the first claim of the thesis. We consider the matrix 

$$\begin{pmatrix}f_{11}&f_{12}\\ f_{12}&f_{22}\end{pmatrix}=\begin{pmatrix}
\sum_{i=2}^df^{(i)}_{11}x_0^{d-i} &\sum_{i=2}^df^{(i)}_{12}x_0^{d-i} \\
\sum_{i=2}^df^{(i)}_{12}x_0^{d-i} &\sum_{i=2}^df^{(i)}_{22}x_0^{d-i} \end{pmatrix}$$

Let $f^{(j)}$ be the first nonzero summand. We get 
$$f_{11}f_{22}-f_{12}^2 = H_{12}(f^{j})x_0^{2d-2j}+\cdots + H_{12}(f^{(d)})$$
where the intermediate terms  involve also some mixed transvectants among the $f^{(i)}$.
By linearity on each row of the determinant, we get some formulas similar to \eqref{eq:puiseux}.
Expanding the determinant, we find the precise formula that is 
$$f_{11}f_{22}-f_{12}^2 = \sum_{i=j}^d\left(\sum_{k_1+k_2=2i}H_{12}(f^{(k_1)},f^{(k_2)})\right)x_0^{2d-2i}.$$

The Hessian of $f^{(j)}(x_1, x_2)$ vanishes, since it corresponds to the first summand obtained for $i=j$. 
Hence,  up to a linear change of coordinates,  the form $f^{(j)}(x_1, x_2)$ can be assumed to be $x_1^j$ (note that $j\ge 2$). Indeed our statement
is invariant by linear change of coordinates involving $x_1, x_2$.

The first claim is that $x_1^j$ divides $f^{(i)}$ for any $i\ge j$.
This is proved by induction on $i$, the starting case $i=j$ has been just granted. For the next case,  $H_{12}(x_1^j,f^{(j+1)})=0$ implies 
$f^{(j+1)}_{22}=0$, so that $f^{(j+1)}=x_1^jn$ for a linear form $n$.
At the next step one has
\begin{equation}H_{12}(x_1^j,f^{(j+2)})+H_{12}(x_1^jn,x_1^jn) = 0\end{equation}
This implies
$j(j-1)x_1^{j-2} f^{(j+2)}_{22} - 2 j^2x_1^{2j-2}n_2^2=0$
hence
$ f^{(j+2)}_{22}$ is proportional to $x_1^j$, which implies that
$ f^{(j+2)}=x_1^jq$ for a quadratic form $q$.
At the next step one has
\begin{equation}H_{12}(x_1^j,f^{(j+3)})+H_{12}(x_1^jn,x_1^jq) = 0\end{equation}
which implies
$$j(j-1)x_1^{j-2} f^{(j+3)}_{22}+\left(j(j-1)x_1^{j-2}n+2jx_1^{j-1}n_1\right)x_1^jq_{22}
-2\left(jx_1^{j-1}q_2+x_1^jq_{12}\right)jx_1^{j-1}n_2=0$$
hence  $f^{(j+3)}_{22}$ is divisible by $x_1^j$, which implies that 
$ f^{(j+3)}=x_1^jc$ for a cubic form $c$. 
Continuing in this way we can prove that $f^{(i)}_{22}$ is divisible by $x_1^j$ for any $i\ge j$,
hence $f^{(i)}$ is divisible by $x_1^j$ (in particular by $x_1^2$) for any $i\ge j$.

Now assume by contradiction that $x_2$ appears with positive exponent in  $f^{(i)}(x_1, x_2)$ for some $i=3,\ldots, d$. Let $x_2^M$ be the maximum appearance, hence $M\ge 1$. It may appear more than once, so that for convenient scalars $a_i$  we have
$$\sum_{i=2}^df^{(i)}x_0^{d-i} = \left(\sum_{i=M+2}^{d}a_ix_0^{d-i}x_1^{i-M}\right)x_2^M+\textrm{lower terms in\ }x_2,$$ 
and $i-M\ge 2$ since $f^{(i)}$ is divisible by  $x_1^2$.
Let $N=\max\{i=M+2,\ldots, d|a_i\neq 0\}$.
Then
$$f_{11}=\sum_{i=M+2}^N (i-M)(i-M-1)a_ix_0^{d-i}x_1^{i-M-2}x_2^M+\textrm{lower terms in\ }x_2,$$
$$f_{12}=\sum_{i=M+2}^N (i-M)Ma_ix_0^{d-i}x_1^{i-M-1}x_2^{M-1}+\textrm{lower terms in\ }x_2,$$
$$f_{22}= \sum_{i=M+2}^N M(M-1)a_ix_0^{d-i}x_1^{i-M}x_2^{M-2}+\textrm{lower terms in\ }x_2$$
so that
$$f_{11}f_{22}-f_{12}^2=-a_N^2M(N-1)(N-M)x_0^{2d-2N}x_1^{2N-2M-2}x_2^{2M-2}+\textrm{lower terms in\ } x_2,$$
where the monomial order is the lexicographical order with $x_2>x_1>x_0$.
The hypothesis implies $a_N=0$ which is the desired contradiction,  that proves the first assertion.

 Now we go on assuming $m=x_1$.  If ${\rm Hess}(f)=0$   then $f=0$  is a cone by  Hesse's Theorem, and it is immediate to check that then $l$ is a multiple of $x_1$. Conversely,
if $l$, $m$ are proportional then  $f=0$ is a cone and  ${\rm Hess}(f)=0$.  

We  finally  show that ${\rm Hess}(f)$ is divisible by $x_0^{2d-4}$. Recall that
$$f=x_0^d+x_0^{d-1}l(x_1,x_2)+\sum_{i= 2}^dc_ix_0^{d-i}x_1^i.$$
If $c_i=0$ for all $2\leq i\leq d$, then $f=0$ is easily seen to be a cone and ${\rm Hess}(f)=0$, in which case the assertion is trivial. If there is an $i$ with $2\leq i\leq d$ such that $c_i\neq 0$, then
we have $f_{11}\neq 0$, 
and also
$$f_{12}=0, \,\, f_{22}=0, \,\, f_{02}=\gamma x_0^{d-2}$$
with $\gamma$ a constant.  Hence the Hessian matrix is

$$\begin{pmatrix}f_{00}&f_{01}&\gamma x_0^{d-2}\\
f_{01}&f_{11}&0\\
\gamma x_0^{d-2}&0&0\end{pmatrix}$$
and the thesis follows. 
\qedd \end{proof}
\vskip 0.6cm
\begin{rem0}\label{rem:maxminor_hessian} {\rm
The plane curves of degree $d$  with equation $f=0$; with $f$ as in  Lemma \ref{lem:heszero}
have a point $P$ of multiplicity $d-1$ (which is $(0:0:1)$ in the homogeneous coordinates of the Lemma if $l=x_2$), which is a hyperflex, in the sense that  all lines
through $P$ meet the curve with multiplicity $d-1$ except the flex tangent (which is $x_0=0$) which meets with the highest multiplicity $d$ at $P$. Note these curves are rational and in some sense are the irreducible curves ``closest'' to the cones. Indeed a plane curve of degree $d$ is a cone if and only if it has a point of multiplicity $d$,  where all derivatives of degree $d-1$ of the defining polynomial vanish.
Instead, for the curves in question,  at the point of multiplicity $d-1$ all derivatives up to order $d-2$ vanish, moreover the tangent cone has degree $d-1$
and consists of a multiple line,  which  means that all derivatives of order $d-1$ vanish except (in the above coordinate system)
$\frac{\partial^{d-1}f}{\partial x_2^{d-1}}(0,0,1)\neq 0$.}

\end{rem0}

\begin{thm0}\label{thm:divides}
Let $r=2$, let $f(t)$ be as in  (\ref{eq:fPuiseux}). Assume $d\geq 4$. Then
 $x_0^{d-3}$ divides the limit of ${\rm Hess}(f(t))$ for $t\to 0$.
\end{thm0}
\begin{proof}
The proof is a case by case analysis of the expansion  (\ref{eq:puiseux}). 

First case, the limit in (\ref{eq:puiseux}) for $t\to 0$ is a summand in the second row.
In this case it is clear that $x_0^{d-2}$ divides the limit.

Second case,  the limit in (\ref{eq:puiseux}) for $t\to 0$ is  a summand in the third row. Let
$H(f_i, f_j, f_k)$ be such a summand.
We have now a few subcases.

\begin{enumerate}
\item{} If $i=j=k$ then $H_{12}(f_i, f_i)$ vanish (since it is a previous summand in  (\ref{eq:puiseux}))
and by Lemma \ref{lem:heszero} we get that $H(f_i, f_i, f_i)$ is divisible by  $x_0^{2d-4}$. 

\item{} If $i<j=k$ then $H_{12}(f_i, f_i) = H_{12}(f_i, f_j) = H_{12}(f_j, f_j) = 0$ and also
$H(f_i, f_i, f_i) = 0$. By  Lemma \ref{lem:hesffg} below we get that $H(f_i, f_j, f_j)$ is divisible by $x_0^{2d-4}$.

\item{} If $i=j<k$ then $H_{12}(f_i, f_i) = H_{12}(f_i, f_j)  = 0$ and also
$H(f_i, f_i, f_i) = 0$.  In this case $H(f_i, f_i, f_k)$ is  divisible by $x_0^{2d-4}$, by Lemma \ref{lem:hesffg} below. 

\item{} If $i<j<k$ then $H_{12}(f_i, f_i) = H_{12}(f_i, f_j) = H_{12}(f_j, f_j) =  H_{12}(f_i, f_k) = H_{12}(f_j, f_k)=  0$ and also
$H(f_i, f_i, f_i) = H(f_i, f_i, f_j) = H(f_i, f_j, f_j) = 0$.
By Lemma \ref{lem:hesfgh} below we get that $H(f_i, f_j, f_k)$ is divisible by  $x_0^{d-3}$.

\end{enumerate}
This concludes the proof of Theorem \ref{thm:divides} after the following two lemmata are proved. \qedd\end{proof}

\begin{lemma0}\label{lem:hesffg}
If $H_{12}(f,f)=0$, $H_{12}(f,g)=0$,  $H(f)=0$,
then $H(f,f,g)$ is divisible by $x_0^{2d-4}$. Moreover, if also   $H_{12}(g,g)$  vanishes then $H(f,g,g)$ is divisible by  $x_0^{2d-4}$. 
\end{lemma0}
\begin{proof} Let us prove the first assertion.
By Lemma \ref{lem:heszero}, the assumptions imply that
$$f=x_0^d+\sum_{i= 1}^dc_ix_0^{d-i}m(x_1,x_2)^i.$$  

The statement is invariant by a linear change of coordinates in $x_1$, $x_2$, hence we may assume $m(x_1,x_2)=x_1$.
It follows that $f_2=0$. We get  $H_{12}(f,g)=\frac{1}{2}f_{11}g_{22}$, hence  either $f_{11}=0$ or $g_{22}=0$.
 
Suppose first $f_{11}=0$. As
$$
f_{11}=\sum_{i= 2}^di(i-1)c_i x_0^{d-i}x_1^{i-2},
$$
we have that $c_i=0$ for all $2\leq i\leq d$, and therefore $f=x_0^d+c_1x_0^{d-1}x_1$, so
$f_{01}$ is divisible by $x_0^{d-2}$. A straightforward calculation shows that $H(f,f,g)=-\frac{1}{3}g_{22}f_{01}^2$, that proves the assertion. 

If $g_{22}=0$, taking into account that $f_2=0$, one easily checks that $H(f,f,g)=0$ and the assertion follows again. 

Let us now prove the second assertion. Again, by Lemma \ref{lem:heszero}  applied to $f$, we have that
\begin{equation}\label{eq:1}
f=x_0^d+\sum_{i= 1}^dc_ix_0^{d-i}(ax_1+bx_2)^i.
\end{equation}
By Lemma \ref{lem:heszero} applied to $g$ we         may assume that
\begin{equation}\label{eq:2} g=x_0^d+x_0^{d-1}l+\sum_{i=2}^d e_ix_0^{d-i}x_1^i
\end{equation}
with $l=l(x_1,x_2)$ a suitable linear form, so that $g_{12}=g_{22}=0$ and $x_0^{d-2}$ divides $g_{02}$. Then again $H_{12}(f,g)=  \frac 12  f_{22}g_{11}=0$, so that either $g_{11}=0$ or $f_{22}=0$.

Assume first that $g_{11}=0$. One has
$$
g_{11}=\sum_{i=2}^d e_ii(i-1)x_0^{d-i}x_1^{i-2},
$$
hence $g_{11}=0$ implies that $e_i=0$ for all $2\leq i\leq d$, so that 
\begin{equation}\label{eq:3}
g=x_0^d+x_0^{d-1}l,
\end{equation}
so $g_{00}$ is divisible by $x_0^{d-3}$ and $g_{01}$ , $g_{02}$ are   divisible by $x_0^{d-2}$. Then the only non--zero summands in $H(f,g,g)$ are as follows
$$3H(f,g,g)=\left|\begin{matrix}g_{00}&g_{01}&g_{02}\\f_{10}&f_{11}&f_{12}\\
g_{20}&0&0\end{matrix}\right|+
\left|\begin{matrix}g_{00}&g_{01}&g_{02}\\
g_{10}&0&0\\
f_{20}&f_{21}&f_{22}\end{matrix}\right|$$ and, taking into account the previous divisibilities,  the assertion follows.

Next we assume $f_{22}=0$. We have
$$
f_{22}=\sum_{i= 2}^dc_i i (i-1) b^2x_0^{d-i} (ax_1+bx_2)^{i-2}.
$$
So, either $b=0$ or $c_i=0$ for all $2\leq i\leq d$. Suppose first that $b=0$. So we may write $f$ as 
\begin{equation}\label{eq:4}
f=x_0^d+\sum_{i= 1}^dc'_ix_0^{d-i}x_1^i
\end{equation} 
so that $f_2=0$.  
Then  $H(f,g,g)$ reduces to 
$$3H(f,g,g)=
\left|\begin{matrix}g_{00}&g_{01}&g_{02}\\
f_{10}&f_{11}&0\\
g_{02}&0&0\end{matrix}\right|$$ and, taking into account that $x_0^{d-2}$ divides $g_{02}$, the assertion follows again.

Finally assume that $c_i=0$ for all $2\leq i\leq d$, so that 
\begin{equation}\label{eq:5} 
f=x_0^d+c_1x_0^{d-1}(ax_1+bx_2),
\end{equation}
and $f_{02}$ is divisible by $x_0^{d-2}$ and $f_{11}=f_{12}=f_{22}=0$. Then one has
$$3H(f,g,g)=\left|\begin{matrix}g_{00}&g_{01}&g_{02}\\g_{10}&g_{11}&0\\
f_{20}&0&0\end{matrix}\right|+
\left|\begin{matrix}f_{00}&f_{01}&f_{02}\\
g_{10}&g_{11}&0\\
g_{20}&0&0\end{matrix}\right|$$
and the assertion again follows. \qedd
\end{proof} 

\begin{lemma0}\label{lem:hesfgh}
If $H_{12}(f,f)=H_{12}(f,g)=H_{12}(g,g)=H_{12}(f,h)=H_{12}(g,h)=0$,  $H(f)=0$,
then $H(f,g,h)$ is divisible by $x_0^{d-3}$. 
\end{lemma0}
\begin{proof} The analysis is similar to the one in the proof of Lemma \ref {lem:hesffg}. 
As in the proof of that lemma, we have that $f$ and $g$ are as in \eqref{eq:1} and \eqref {eq:2} respectively, so that $g_{12}=g_{22}=0$ and $x_0^{d-2}$ divides $g_{02}$. Again $H_{12}(f,g)=f_{22}g_{11}=0$, so that either $g_{11}=0$ or $f_{22}=0$.

Assume first that $g_{11}=0$. As in the proof of Lemma \ref {lem:hesffg}, this yields that $g$ is as in \eqref {eq:3},
so $g_{00}$ is divisible by $x_0^{d-3}$ and $g_{01}$   is divisible by $x_0^{d-2}$. Now let us look at the six determinants that appear as summands in $H(f,g,h)$. In the first and the last ones the row $(g_{10},g_{11},g_{12})=(g_{10},0,0)$ appears and $g_{10}$ is divisible by $x_0^{d-2}$. In the second and the fifth ones the row $(g_{20},g_{21},g_{22})=(g_{20},0,0)$ appears and $g_{20}$ is divisible by $x_0^{d-2}$. In the third and the fourth ones the row  $(g_{00},g_{01},g_{02})$ appears and all three entries are divisible by  $x_0^{d-3}$. This proves that $H(f,g,h)$ is divisible by  $x_0^{d-3}$.  

Assume next that $f_{22}=0$. As in the proof of Lemma \ref {lem:hesffg}, this yields that 
either $b=0$ or $c_i=0$ for all $2\leq i\leq d$. 

Suppose first that $b=0$. Again as in the proof of Lemma \ref {lem:hesffg}, this implies that 
$f$ can be written as in \eqref {eq:4}, 
so that $f_2=0$.  Then we have $H_{12}(f,h)=f_{11}h_{22}=0$  and $H_{12}(g,h)=g_{11}h_{22}=0$. Suppose $h_{22}=0$. We compute
{\tiny \begin{align}\nonumber 6H(f,g,h)=&\\ \nonumber\left|\begin{matrix}f_{00}&f_{01}&0\\
g_{10}&g_{11}&0\\
h_{20}&h_{21}&0\end{matrix}\right|+&
\left|\begin{matrix}f_{00}&f_{01}&0\\
h_{10}&h_{11}&h_{12}\\
g_{20}&0&0\end{matrix}\right|+\left|\begin{matrix}g_{00}&g_{01}&0\\
f_{10}&f_{11}&0\\
h_{20}&h_{21}&0\end{matrix}\right|+
\left|\begin{matrix}g_{00}&g_{01}&0\\
h_{10}&h_{11}&h_{12}\\
0&0&0\end{matrix}\right|+
\left|\begin{matrix}h_{00}&h_{01}&h_{02}\\
f_{10}&f_{11}&0\\
g_{20}&0&0\end{matrix}\right|+\left|\begin{matrix}h_{00}&h_{01}&h_{02}\\
g_{10}&g_{11}&0\\
0&0&0\end{matrix}\right|\end{align}}
so that the only two non--zero summands are divisible by $g_{20}$ that in turn is divisible by $x_0^{d-2}$, and the assertion follows.

If $h_{22}\neq 0$, then we have $f_{11}=g_{11}=0$, which implies that
$$
f=x_0^d+c_1'x_0^{d-1}x_1, \,\,\,\, g=x_0^d+x_0^{d-1}l
$$
and this yieds that $x_0^{d-3}$ divides $f_{00}$ and $g_{00}$ and $x_0^{d-2}$ divides $f_{01}$ and $g_{01}$. Then again we compute
{\tiny \begin{align}\nonumber 6H(f,g,h)=&\\ \nonumber\left|\begin{matrix}f_{00}&f_{01}&0\\
g_{10}&0&0\\
h_{20}&h_{21}&h_{22}\end{matrix}\right|+&
\left|\begin{matrix}f_{00}&f_{01}&0\\
h_{10}&h_{11}&h_{12}\\
g_{20}&0&0\end{matrix}\right|+\left|\begin{matrix}g_{00}&g_{01}&0\\
f_{10}&0&0\\
h_{20}&h_{21}&h_{22}\end{matrix}\right|+
\left|\begin{matrix}g_{00}&g_{01}&0\\
h_{10}&h_{11}&h_{12}\\
0&0&0\end{matrix}\right|+
\left|\begin{matrix}h_{00}&h_{01}&h_{02}\\
f_{10}&0&0\\
g_{20}&0&0\end{matrix}\right|+\left|\begin{matrix}h_{00}&h_{01}&h_{02}\\
g_{10}&0&0\\
0&0&0\end{matrix}\right|\end{align}}
and we see that the non--zero summands are divisible by $x_0^{2d-4}$ and the assertion follows.

Finally we have to analyse the case in which $c_i=0$ for all $2\leq i\leq d$, so that $f$ is as in \eqref {eq:5}, 
that yields $f_{11}=f_{22}=f_{12}=0$ and $x_0^{d-3}$ divides $f_{00}$ whereas $x_0^{d-2}$ divides $f_{01}$ and $f_{02}$. We have again $H_{12}(g,h)=g_{11}h_{22}=0$, so that either $g_{11}=0$ of $h_{22}=0$.

Suppose $g_{11}=0$.  Then, as above, we have $g=x_0^d+x_0^{d-1}l$ so that $x_0^{d-3}$ divides $g_{00}$, $x_0^{d-2}$ divides $g_{01}$ and $g_{02}$. We have 
{\tiny \begin{align}\nonumber 6H(f,g,h)=&\\ \nonumber\left|\begin{matrix}f_{00}&f_{01}&f_{02}\\
g_{10}&0&0\\
h_{20}&h_{21}&h_{22}\end{matrix}\right|+&
\left|\begin{matrix}f_{00}&f_{01}&f_{02}\\
h_{10}&h_{11}&h_{12}\\
g_{20}&0&0\end{matrix}\right|+\left|\begin{matrix}g_{00}&g_{01}&g_{02}\\
f_{10}&0&0\\
h_{20}&h_{21}&h_{22}\end{matrix}\right|+
\left|\begin{matrix}g_{00}&g_{01}&g_{02}\\
h_{10}&h_{11}&h_{12}\\
0&0&0\end{matrix}\right|+
\left|\begin{matrix}h_{00}&h_{01}&h_{02}\\
f_{10}&0&0\\
g_{20}&0&0\end{matrix}\right|+\left|\begin{matrix}h_{00}&h_{01}&h_{02}\\
g_{10}&0&0\\
0&0&0\end{matrix}\right|\end{align}} and again $H(f,g,h)$ is divisible by $x_0^{d-2}$.

Finally, if $h_{22}=0$, we have
{\tiny \begin{align}\nonumber 6H(f,g,h)=&\\ \nonumber\left|\begin{matrix}f_{00}&f_{01}&f_{02}\\
g_{10}&g_{11}&0\\
h_{20}&h_{21}&0\end{matrix}\right|+&
\left|\begin{matrix}f_{00}&f_{01}&f_{02}\\
h_{10}&h_{11}&h_{12}\\
g_{20}&0&0\end{matrix}\right|+\left|\begin{matrix}g_{00}&g_{01}&g_{02}\\
f_{10}&0&0\\
h_{20}&h_{21}&h_{22}\end{matrix}\right|+
\left|\begin{matrix}g_{00}&g_{01}&g_{02}\\
h_{10}&h_{11}&h_{12}\\
f_{20}&0&0\end{matrix}\right|+
\left|\begin{matrix}h_{00}&h_{01}&h_{02}\\
f_{10}&0&0\\
g_{20}&0&0\end{matrix}\right|+\left|\begin{matrix}h_{00}&h_{01}&h_{02}\\
g_{10}&g_{11}&0\\
f_{20}&0&0\end{matrix}\right|\end{align}}
and by the divisibilities noted above, we again have that $H(f,g,h)$ is divisible by $x_0^{d-2}$.\qedd
\end{proof} 
\vskip 0.5cm

\begin{coro0}\label{coro:indet}
\begin{enumerate}

\item{} The point $(x_0^{2k},q^{3(k-1)})$ does not belong to the closure $G_{2k,2}$,
of the graph of the Hessian map $h_{2k,2}$, as in  (\ref{eq:notin}).

\item{} The point $(x_0^{2k+1},q^{3(k-1)}x_0^3)$ does not belong to the closure $G_{2k+1,2}$,
of the graph of the Hessian map $h_{2k+1,2}$, as in  (\ref{eq:notinodd}) if $k\ge 3$.

\item{} The point $(x_0^{2k},q^{3(k-2)}x_0^6)$ does not belong to the closure $G_{2k,2}$,
of the graph of the Hessian map $h_{2k,2}$, as in  (\ref{eq:notineven2}) if $k\ge 5$.
\end{enumerate}
\end{coro0}
\begin{proof}
The first case is clear by Theorem \ref{thm:divides}. When $d=2k+1$, in order to apply Theorem \ref{thm:divides} we need $2k+1-3>3$, so $k\ge 3$.
In the last case we need $2k-3>6$, so $k\geq 5$. \qedd
\end{proof}

\begin{thm0}\label{thm:main} The Hessian map $h_{d,2}$
is birational  onto  its image for plane curves of even degree $d\ge 4$, $d\neq 5$.
\end{thm0}

\begin{proof} We apply Theorem \ref{thm:criterion}.   Let first $d$ be even. The numerical assumption in (\ref{eq:notin}) for $r=2$ is $2m^2+m-3k\neq 0$  for every $m$ such that 
$1\le m\le k$. The equation $2m^2+m-3k=0$ has no integer solutions for $k=2,\ldots, 6$, so in these cases we get the thesis by Theorem \ref{thm:criterion} ,
since the assumption (\ref{eq:notin}) about the graph  is satisfied thanks to
Corollary \ref{coro:indet}, (1).  For $k\ge 7$  the numerical assumption  (\ref{eq:notineven2}) of  Theorem \ref{thm:criterion} for $r=2$ is satisfied thanks to Theorem \ref{thm:integralPts2} of the Appendix (read there $(k,m)=(x,y)$). The assumption  (\ref{eq:notineven2}) about the graph   is satisfied thanks to
Corollary \ref{coro:indet}, (3)  and we get the thesis by  Theorem \ref{thm:criterion}.
Let now $d$ be odd.  The numerical assumption  (\ref{eq:notinodd}) of  Theorem \ref{thm:criterion} for $r=2$ is satisfied thanks to Theorem \ref{thm:integralPts} of the Appendix. The assumption  (\ref{eq:notinodd}) about the graph  is satisfied thanks to
Corollary \ref{coro:indet}, (2),  and we get again the thesis by Theorem \ref{thm:criterion}. 
\qedd
\end{proof}

\section{Appendix: Integral points on two elliptic curves}
\smallskip
\begin{center}by \textsc{Jerson Caro and Juanita Duque-Rosero}\end{center}
\medskip

In this note, we answer a question from Ciro Ciliberto and Giorgio Ottaviani and show the following theorems.

\begin{thm0}\label{thm:integralPts}
    The integral points on the curve $y^2(2x + 1) + y(x + 2) - 3x(x + 1) = 0$ are
    $$\Omega_1\colonequals \{(0, 0), (-1, 1), (1, -2),(1, 1), (-1, 0), (0, -2)\}.$$
\end{thm0}

\begin{thm0}\label{thm:integralPts2}
    The integral points on the curve $2xy^2 + y(-3x + 3) - x(3x - 1) = 0$ are
    $$\Omega_2\colonequals \{(1, -1), (9 , -3), (2 , 2), (1 ,1), (0 , 0), (-1 ,2), (-1 , 1)\}.$$
\end{thm0}

\subsection{The curve  $y^2(2x + 1) + y(x + 2) - 3x(x + 1) = 0$  }\label{sec:firstCurve}
We first change coordinates  using the map $(x:y:z)\mapsto(x:x+y:z)$ in homogeneous coordinates.
To set up some notation, consider the projective curves
$$C\colon 2x^3 + 4x^2y + 2xy^2 - x^2z + 3xyz + y^2z - xz^2 + 2yz^2=0,$$
\begin{equation}\label{EC 1}
W\colon y^2z=x^3 - \frac{35}{16}x^2z + \frac{21}{16}xz^2 + \frac{9}{64}z^3.
\end{equation}
We have that $W$ is isomorphic to the elliptic curve with LFMDB label 
\textrm{366.b1} (see \cite {LMFDB}). The Mordell--Weil rank of $W$ is 1, and the torsion of the Mordell--Weil group is $\Z/3\Z$.  The curves $C$ and $W$ are birational:
$$\begin{array}{ccccc}\rho_1\colon &C&\to &W\\
&(x:y:z)&\mapsto&\displaystyle\left(6xy:
3x^2 - \frac{9}{2}xy - 3y^2 + \frac{3}{2}xz - 6yz: -8x^2 - 4xz\right)
\end{array}$$
A simple computation shows that the only point $P\in C(\Z)$ for which  $\rho_1$ is not defined is $(0:0:1)$.

An integral model for $W$ is the curve $$X\colon y^2z = x^3 - 8960x^2z + 22020096xz^2 + 9663676416z^3.$$
We also have
$$\begin{array}{ccccc}\rho_2\colon &W&\to  &X\\
    &(x : y : 1) &\mapsto &\displaystyle \left(2^{12}x : 2^{18}y : 1\right).
\end{array}$$

Let $C_0$, $W_0$, and $X_0$ be the affine charts where $z=1$.  Then Theorem~\ref{thm:integralPts} can be restated as $C_0(\Z)=\Omega_1$.  We now prove three lemmas that will allow us to identify $C_0(\Z)$ using points on $W$ and $X$.

\begin{lemma0}\label{lem:pointsInC}
    Let $P\colonequals (x,y)\in C_0(\Z)\setminus\{(0,0)\}$.  Then, $\rho_1(P)\in W_0\left(\Z\left[\frac{1}{6}\right]\right)\cup\{(0:1:0)\}$.
\end{lemma0}
\begin{proof}
    Let $(x:y:1)\in C_0(\Z)$. Since $\rho_1$ is defined over $\Z\left[\frac{1}{2}\right]$, then $\rho_1(P)$ may have powers of $2$ as denominators.  The only part that is left is to show that if $-8x^2 - 4x\ne 0$, then normalizing $\rho_1$ to have the coordinate $z$ equal to one produces coordinates $x$ and $y$ in $\Z\left[\frac{1}{6}\right]$.  To prove this it is enough to show that for any prime number $p>3$ and $x,y$ as before,
    $$
    v_p(6xy),v_p\left(
3x^2 - \frac{9}{2}xy - 3y^2 + \frac{3}{2}x - 6y\right) \geq v_p(-8x^2 - 4x).
    $$
    Let $p>3$ be a prime number.  Then
    $$v_p(-8x^2 - 4xz)=v_p(-4x(2x+1))=v_p(x)+v_p(2x+1).$$
    We note that since $p>3$, only one of those valuations can be nonzero.  We first assume that $v_p(x)>0$.  It follows that $v(6xy)\ge v_p(x)$.  For the other coordinate:
    \begin{align*}
        v_p\left(3x^2 - \frac{9}{2}xy - 3y^2 + \frac{3}{2}x - 6y\right)&=v_p\left(2x^2-3xy-2y^2+x-4y\right)\\
        &= v_p\left(3xy+4x^3+8x^2y+4xy^2-x\right)\\
        &=v_p(x)+v_p(3y+4x^2+8xy+4y^2-1)\\
        &\ge v_p(x),
    \end{align*}
    where the second equality follows from the equation of $C$.

    Now we assume that $v_p(2x+1)>0$ and we note that this implies that $v_p(x)=0$.
    By the equation defining $W$, we have that modulo $p$ the only point at infinity is $[0:1:0]$, in particular, $v_p(y)=v_p(6xy)\geq v_p(-8x^2 - 4x)=v_p(2x+1)$. On the other hand, we have:
    {\footnotesize
    \begin{align*}
        v_p\left(3x^2 - \frac{9}{2}xy - 3y^2 + \frac{3}{2}x - 6y\right)&=v_p\left((2x+1)\left(\frac{3}{2}x\right)-y\left(\frac{9}{2}x+3y+6\right)\right)\\
        &\ge \min\left\{v_p(2x+1)+v_p\left(\frac{3}{2}x\right),v_p(y)+v_p\left(\frac{9}{2}x+3y+6\right)\right\}\\
        &\ge\min\{v_p(2x+1),v_p(y)\}\\
        &\ge v_p(2x+1),
    \end{align*}}
    which yields the desired result.\qed
\end{proof}

\begin{lemma0} \label{lem:pointsInW}
    Let $x,y\in \Z\left[\frac{1}{6}\right]$ and let $P\colonequals(x,y)$.  Then $P\in W_0\left(\Z\left[\frac{1}{6}\right]\right)$ if and only if $\rho_2(P)\in X_0\left(\Z\left[\frac{1}{6}\right]\right)$.
\end{lemma0}
\begin{proof}
    The map $\rho_1$ is already normalized, and the only denominators are powers of $2$. \qed
\end{proof}

\begin{lemma0}\label{lem:ptsW}
    For $W_0$ as above, we have \begin{align*}
    \left.W_0\left(\Z\left[\frac{1}{6}\right]\right)=\right\{&\left(-\frac{29}{324} , \pm\frac{817}{11664} \right), \left(0 , \pm\frac{3}{8}\right), \left(\frac{3}{16} , \pm\frac{9}{16} \right), \left(\frac{3}{4} , \pm\frac{9}{16} \right),\\
    &\left(\frac{5}{4} , \pm\frac{9}{16} \right), \left(\frac{3}{2} , \pm\frac{3}{4} \right), \left(\frac{2145}{1024} , \pm\frac{51633}{32768} \right), \left(3 , \pm\frac{27}{8} \right), \\&\left.\left(\frac{21}{4} , \pm\frac{153}{16} \right), \left(27 , \pm\frac{1077}{8} \right)\right\}.
\end{align*}
\end{lemma0}
\begin{proof}
    Using the \magma~function \verb|SIntegralPoints| \cite{Magma}, we compute the set  $X_0\left(\Z\left[\frac{1}{6}\right]\right)$ and then take the preimage under $\rho_2$ to $W$.  By Lemma \ref{lem:pointsInW}, this set equals $W_0\left(\Z\left[\frac{1}{6}\right]\right)$.\qed
\end{proof}\medskip

\begin{proof}[Proof of Theorem~\ref{thm:integralPts}]
We first compute integral points on $C_0$.
Lemma~\ref{lem:pointsInC} implies that the images via $\rho_1$ of  the points in $C_0(\Z)\setminus\{(0,0)\}$ lie in
$$
W_0\left(\Z\left[\frac{1}{6}\right]\right)\cup \{(0:1:0)\}.
$$
The only points in $C_0(\Z)$ with coordinates $x$ or $y$ equal to $0$ are $(0,0)$, $(0,-2)$, and $(1,0)$.  Consequently, by the equations defining $\rho_1$, the other integral points of $C_0(\Z)$ map to $W_0(\Z\left[\frac{1}{6}\right])$ via $\rho_1$.

For $(a,b)\in W_0(\Z\left[\frac{1}{6}\right])$ and $x,y\in \overline{\QQ}^\times$, we have that $\rho_1(x,y)=(a,b)$ if and only if
\[
a=\frac{-3y}{4x+2}\quad \text{and} \quad b=\frac{3x^2 - \frac{9}{2}xy - 3y^2 + \frac{3}{2}x - 6y}{-8x^2-4x}.
\]
We solve for $y$ obtaining $\frac{-a(4x+2)}{3}$ and the possible values of $x$ are the roots of the quadratic polynomial
\[
x^2\left(3+6a-\frac{16a^2}{3}+8b\right)
+x\left(11a-4b-\frac{16a^2}{3}+\frac{3}{2}\right)
+4a-\frac{4a^2}{3}.
\]
Using these identities, we find the values for $x$ and $y$ associated with all points obtained from Lemma~\ref{lem:ptsW}.  Then we check that the only integer values are:
$$\{(0, 0), (-1, 0), (1, -3),(1, 0), (-1, 1), (0, -2)\}.$$
We recall that the morphism from the curve of Theorem~\ref{thm:integralPts} to $C$ is $(x:y:z)\mapsto(x:x+y:z)$.  That allows us to recover the desired set $\Omega_1$.
\end{proof}

\subsection{The curve $2xy^2 + y(-3x + 3) - x(3x - 1) = 0$}

We follow the same ideas as in Section~\ref{sec:firstCurve}.  We consider the curve $$C\colon 2xy^2 + y(-3x + 3) - x(3x - 1) = 0,$$
with Weierstrass model
\begin{equation}\label{EC 2}
W\colon y^2 = x^3 + \frac{1}{4}x^2 - 27x + 81.    
\end{equation}
The curve $W$ is isomorphic to the elliptic curve with LFMDB label 
\textrm{1002.e1}.  The Mordell-Weil rank is 1, and this group has no torsion.  The curves $C$ and $W$ are birational:
$$\begin{array}{ccccc}\rho_1\colon &C&\to &W\\
&(x:y:z)&\mapsto&\displaystyle\left(6yz: 12y^2 - 9xz -9yz:-xz\right)
\end{array}$$

An integral model for $W$ is the curve $$X\colon y^2 = x^3 + 4x^2 - 6912x + 331776.$$
We also have
$$\begin{array}{ccccc}\rho_2\colon &W&\to  &X\\
    &(x : y : 1) &\mapsto &\displaystyle \left(16x : 64y : 1\right).
\end{array}$$

Let $C_0$, $W_0$, and $X_0$ be the affine charts where $z=1$.  Then Theorem~\ref{thm:integralPts2} can be restated as $C_0(\Z)=\Omega_2$.  We now prove lemmas similar to the ones in Section~\ref{sec:firstCurve} to identify $C_0(\Z)$ using points on $W$ and $X$.

\begin{lemma0}\label{lem:pointsInC2}
    For $P\colonequals(x,y)\in C_0(\Z)$ with $x\ne 0$, we have $\rho_1(P)\in W_0\left(\Z\right)$.
\end{lemma0}
\begin{proof}
    Let $P\colonequals(x:y:1)\in C_0(\Z)$ with $x\ne 1$.  Since $\rho_1$ is defined over $\Z$ and $-xz\ne 0$, we need to show that normalizing $\rho_1(P)$ to have the coordinate $z$ equal to 1 produces coordinates $x$ and $y$ in $\Z$.  The fact we want follows from the equality:
    $$v_p(3y)= v_p(x)+v_p(2y^2 - 3x - 3y + 1 ).$$
    Even when $p=3$, this is enough since $3|6x$ and $3|(12y^2-9x-9z)$. \qed
\end{proof}

\begin{lemma0} \label{lem:pointsInW2}
    If $P\colonequals(x,y)\in W_0\left(\Z\right)$  then $\rho_2(P)\in X_0\left(\Z\right)$.
\end{lemma0}
\begin{proof}
    Clear from the definition of $\rho_2$.\qed 
\end{proof}

\begin{lemma0}\label{lem:ptsW2}
    For $W_0$ as above, we have 
    {\footnotesize \begin{align*}  W_0(\Z)=\{(-6 : \pm6 : 1), (0 : \pm9 : 1), (2 : \pm6 : 1), (6 : \pm12 : 1), (12 : \pm39 : 1), (54 : \pm396 : 1)\}.
\end{align*}}
\end{lemma0}
\begin{proof}
    Using the \magma~function \verb|IntegralPoints| \cite{Magma}, we compute the set  $X_0\left(\Z\right)$ and then take the preimage under $\rho_2$ to $W$.  By Lemma \ref{lem:pointsInW}, this set contains $W_0\left(\Z\right)$, so we pick the points with integral coordinates. \qed
\end{proof}\medskip

\begin{proof}[proof of Theorem~\ref{thm:integralPts2}]
By Lemma~\ref{lem:pointsInC2}, we have that $\rho_1(C_0(\Z)\setminus\{(0,0)\})\subseteq W_0\left(\Z\right).$
For $(a,b)\in W_0(\Z)$ and $x,y\in \overline{\QQ}$ the equality $\rho_1(x,y)=(a,b)$ holds if and only if
\[
a=\frac{-6y}{x}\quad \text{and} \quad b=\frac{9x+9y-12y^2}{x}.
\]
We note that $x\ne 0$ since this would also imply that $y=0$, and the image under $\rho_1$ is not defined.
Using the identities above, we solve for $(x,y)$ from each value of $(a,b)$ in Lemma~\ref{lem:ptsW}.  Then we add $(0,0)$ since we excluded it in Lemma~\ref{lem:pointsInC2}.  Finally, we check that the only integer values are the points in $\Omega_2$. \qed
\end{proof}

{\begin{rem0}\label{rem:function_magma} {\rm
The \magma~functions \verb|IntegralPoints| and \verb|SIntegralPoints| \cite{Magma} implement a deterministic method based on the work of Stroeker and Tzanakis \cite{ST}, which uses elliptic logarithms. This method is deterministic provided the group structure of the rational points on the elliptic curve is fully determined, which, as mentioned above, is indeed the case here. Specifically, the group of rational points on the elliptic curve \eqref{EC 1} is isomorphic to $\mathbb{Z} \oplus \mathbb{Z}/3\mathbb{Z}$, while the group of rational points on the elliptic curve \eqref{EC 2} is isomorphic to $\mathbb{Z}$.
}
\end{rem0}}


\bigskip
{\small


\begin{thebibliography} {Dillo}
\bibitem[Beo]{Beo} V. Beorchia, {\em Generic injectivity of Hessian maps of ternary forms},
arXiv:2406.05423.
\bibitem[BK]{BK} A. I. Bondal, M.M. Kapranov, {\em Homogeneous bundles} in: Helices and Vector
Bundles, London Math. Soc. Lecture Note Ser. 148, Cambridge Univ. Press,
Cambridge, 1990, 45 -- 55.
\bibitem[BFP]{BFP} D. Bricalli, F. Favale, G. P. Pirola, {\em A theorem of Gordan and Noether via Gorenstein rings},
Selecta Math. (N.S.) 29 (2023), no. 5, Paper No. 74, 25 pp.
\bibitem[CO]{CO} C. Ciliberto, G. Ottaviani, {\em The Hessian map}, Int. Math. Res. Notices, 8 (2022), 5781--5817.
\bibitem[CRS]{CRS} C. Ciliberto, F. Russo, A. Simis, {\em  Homaloidal hypersurfaces and hypersurfaces with vanishing Hessian},
Adv. Math. 218 (2008), no.6, 1759--1805.
\bibitem[CNOS]{CNOS} A. Conca, S. Naldi, G. Ottaviani, B. Sturmfels, {\em Taylor Polynomials of Rational Functions}, Acta Math. Vietnamica, 
 49 (2024), no. 1, 19--37.
\bibitem[Dol]{Dol} I. Dolgachev, {\em Polar Cremona transformations},
Michigan Math. J.48 (2000), 191--202.

\bibitem[Dolg]{Dolg} I. Dolgachev, {\em Classical Algebraic Geometry: a modern view}, Cambridge University Press, 2011.

\bibitem[FH]{FH} W. Fulton and J. Harris, {\em  Representation theory}, Graduate Texts in Mathematics 129. Springer-
Verlag, 1991.

 \bibitem[GW]{GW} R. Goodman, N. Wallach, {\em Representations and invariants of the classical groups},
Encyclopedia Math. Appl., 68, Cambridge, 1998.

\bibitem[Kap]{Kap} M. M. Kapranov, {\em On the derived categories of coherent sheaves on some homogeneous spaces},
Invent. Math. 92 (1988), no. 3, 479--508.

\bibitem[KT]{KT} K. Koike, I. Terada, {\em Young-diagrammatic methods for the representation theory of the classical groups of type  $B_n$, $C_n$, $D_n$}, 
J. Algebra 107 (1987), no. 2, 466--511.

\bibitem[Lit]{Lit} P. Littelmann, {\em  A generalization of the Littlewood-Richardson rule},
 J. Algebra 130 (1990), no. 2, 328--368.

\bibitem[LMFDB]{LMFDB} The LMFDB Collaboration, The L-functions and modular forms database, https://www.lmfdb.org, 2024, 


\bibitem[Magma]{Magma}
Wieb Bosma, John Cannon, and Catherine Playoust, \emph{The Magma algebra system.\ I.\ The user language}, J.\ Symbolic Comput.\ \textbf{24} (1997), vol.\ 3--4, 235--265.

\bibitem[Man]{Man} L. Manivel, {\em Prehomogeneous spaces and projective geometry}, Rend. Semin. Mat. Univ. Politec. Torino71 (2013), no.1, 35--118.

\bibitem[O12]{O12} G. Ottaviani, {\em  Five Lectures on projective Invariants}, lecture notes for Trento school, September 2012, Rendiconti del Seminario Matematico Univ. Politec. Torino, vol. 71, 1 (2013), 119--194.

\bibitem[OR]{OR} G. Ottaviani, E. Rubei, {\em Quivers and the cohomology of homogeneous vector bundles},
Duke Math. J. 132 (2006), no. 3, 459--508.

\bibitem[Ru]{Ru}  F.~Russo,
  {\em On the geometry of some special projective varieties},
    Lecture Notes of the Unione Matematica Italiana, {\bf 18}, Springer, Cham, 2016.

\bibitem[Seg]{Seg} B. Segre, {\em  Bertini forms and Hessian matrices},
J. London Math. Soc. 26 (1951), 164--176.

\bibitem[ST]{ST} R. Stroeker, N. Tzanakis (1994). \emph{Solving elliptic diophantine equations by estimating linear forms in elliptic logarithms.} Acta Arithmetica, 67(2), 177-196.

\end{thebibliography}}
\end{document}